\documentclass[preprint]{imsart}

\usepackage[numbers]{natbib}
\usepackage[OT1]{fontenc}
\usepackage{amssymb, amsmath, amsfonts,amsthm}
\usepackage[colorlinks]{hyperref}
\usepackage{graphicx}
\usepackage{amsfonts}
\usepackage{epstopdf}
\usepackage{subfigure}
\usepackage{amssymb}
\usepackage{amsmath}
\usepackage{url}
\usepackage{accents}
\usepackage{rotating,color}
\usepackage[dvipsnames, svgnames]{xcolor}
\usepackage{tikz}    

\arxiv{math.PR/2502.09185.v2}

\startlocaldefs
\numberwithin{equation}{section}
\theoremstyle{plain}
\newtheorem{theorem}{Theorem}[section]

\newtheorem{corollary}{Corollary}[section]
\newtheorem{lemma}{Lemma}[section]

\newtheorem{remark}{Remark}[section]
\endlocaldefs

\font\eufm = eufm10 scaled 1000

\def\num#1{\hbox{(\ref{#1})}}
\def\d{\displaystyle }

\def\Be{\begin{equation}}
\def\Eeq{\end{equation}}
\def\be{\begin{eqnarray}}
\def\eeq{\end{eqnarray}}
\def\ben{\begin{enumerate}}
\def\een{\end{enumerate}}
\def\ba{\begin{array}}
\def\ea{\end{array}}
\def\bt{\begin{tabular}}
\def\et{\end{tabular}}
\def\bc{\begin{center}}
\def\ec{\end{center}}
\def\bi{\begin{itemize}}
\def\ei{\end{itemize}}
\def\bd{\begin{document}}
\def\ed{\end{document}}
\def\nn{\nonumber\\ }
\def\nnn{\nonumber\eeq }
\def\lf{\lefteqn}
\DeclareMathOperator*{\argmax}{arg\,max}
\DeclareMathOperator*{\argmin}{arg\,min}
\def\choose#1#2{{\;\big( \ba{c} \vspace*{-1mm} #1 \\ \vspace*{1mm} #2 \ea\big)\;}}
\def\as{\,\mbox{ as }\,}
\def\If{\,\mbox{ if }\,}
\def\for{\,\mbox{ for }\,}
\def\and{\,\mbox{ and }\,}
\def\And{\,\mbox{ and }\,}
\def\textas{\,\mbox{ as }\,}
\def\textif{\,\mbox{ if }\,}
\def\textfor{\,\mbox{ for }\,}
\def\textand{\,\mbox{ and }\,}
\def\P{\mathbb P}
\def\zero{0\!\!\!\;\!\! 0}
\def\X{\mbox{\boldmath $X$}}
\def\Y{\mbox{\boldmath $Y$}}
\def\D{\mathbb{V\textnormal{ar}}}
\def\E{\mathbb E}
\def\R{\mathbb R}
\def\N{\mathbb N}
\def\Z{\mathcal Z}
\def\var{{\,{\mathrm{Var}}}}
\def\VAR{{\,\mathrm{VAR}}}
\def\cov{{\,{\mathrm{Cov}}}}
\def\std{{\,\mathrm{Std}}}
\def\corr{{\,\mathrm{Corr}}}
\def\x{\mbox{\boldmath $x$}}
\def\ep{\varepsilon}
\def\til{\!\!\!_{_\sim}}
\def\haha{$\stackrel{..}{\smile}$}
\def\edk{\hbox{\rule{2mm}{2mm}}\medskip\smallskip}
\def\l{\left\{ }
\def\r{\right\} }
\def\ll{\left.}
\def\rr{\right.}
\def\ra{\rightarrow}
\def\supp{\mbox{supp}}
\def\c{ \left|\ba{c} \vspace*{-5mm}\\ \ea\right.\!\!\!\!\!\! }
\def\cc{ \,\left|\ba{c} \vspace*{-5mm}\\ \ea\right.\!\!\!\!\! }
\def\*{$\cdot$}
\def\Xbar{\bar{X}}
\def\topw#1#2{\xrightarrow[\!\!#2\!\!]{#1}}

\begin{document}
\def\KL{\mbox{\rm KL}}
\def\MA{\mbox{\eufm{M}}}
\def\N{\mathbb{N}}
\def\Deltab{{\boldsymbol \Delta}}
\def\di{\displaystyle}
\def\indi{1\hspace{-1,1mm}{\rm I}}
\def\Fb{{\mathbb{F}}}
\def\Ff{\mbox{\eufm F}}
\def\Rb{{\mathbb{R}}}
\def\TA{\mbox{\eufm{T}}}
\def\Tb{{\mathbb{T}}}
\def\Pr{{\mathbb{P}}}
\def\F{{\cal{F}}}
\def\zetab{\mbox{\boldmath$\zeta$}}
\def\Pc{{\cal{P}}}
\def\E{{\mathbb E}}
\def\Var{{\mathrm {Var}}}
\def\alph{a}
\def\bet{b}
\def\X{\mbox{\boldmath $X$}}
\def\Ybi{\mbox{\textbf{\textit Y\rule{-.6ex}{0em}}}}
\def\Nb{{\mathbb N}}
\def\ubi{\mbox{\textbf{\textit u}}}
\def\ubis{\mbox{\textbf{\scriptsize\textit u}}}
\def\Zb{{\mathbb Z}}

\def\EM{{\cal E}}   
\def\M{{E}}   
\def\V{{V}}   
\def\lam{\lambda}   
\def\lamstar{\lam^{\!\ast}}
\def\ttau{\tau^{\circ}}
\def\MTBFA{\mbox{MT2FA}}
\def\MD{\mbox{MD}}
\def\ARL{\mbox{ARL}}

\arraycolsep=2pt


\begin{frontmatter}
\title{{\large\bf
PRACTICAL PROPERTIES OF THE CUSUM PROCESS
}}
\runtitle{Practical Properties of the CUSUM process}

\begin{aug}
\author{\fnms {Michael} \snm{Baron$^{1,}$}\thanksref{t1}\ead[label=e1]{baron@american.edu}}
\author{\fnms{Sergey V.} \snm{Malov$^{2,3,}$}\thanksref{t2}\ead[label=e2]{sergey.v.malov@gmail.com }}

\thankstext{t2}{E-mail: sergey.v.malov@gmail.com }
\thankstext{t1}{E-mail: baron@american.edu }

\runauthor{ M. Baron  $\&$ S.V. Malov }

\address{$^1$American University, Washington D.C., USA}
\address{$^2$Peter the Great St.-Petersburg Polytechnic University, St.~Petersburg, Russia}
\address{$^3$St.-Petersburg Electrotechnical University, St.~Petersburg, Russia}

\end{aug}

\begin{abstract}
We establish new properties of the Lindley process, also known as the waiting time process or the cumulative-sum (CUSUM) process, as well as its running maximum. Unlike many other authors,
the study focuses on exact and asymptotic expressions for the moment generating function (MGF) of CUSUM. It includes precise finite-sample expressions for the MGF and moments of the process, along with fast recursive computing algorithms, lower and upper bounds, and the classification of large-sample asymptotics of the CUSUM MGF. Results are applied to single, multiple, and transient change-point problems, for the calculation of thresholds that provide a desired control of familywise false alarm rates, as well as the quantiles of queuing processes and probabilities of their large deviation at least once over a given finite time interval.
 \end{abstract}



\end{frontmatter}


\section{Introduction}\label{sec:introduction}
A cumulative-sum (CUSUM) process, based on a sequence of random variables $\{ Y_n, n=1,2,\ldots \}$, is defined recursively as
  \Be\label{CUSUMdefinition}
  W_0=0, \ W_{n+1} = \max\{ W_n + Y_{n+1}, 0 \}.
  \Eeq
This construction appears in multiple important areas of applications such as queuing theory, quality and process control, and
change-point detection.

Kendall \cite{Kendall53}, Lindley \cite{Lindley1952queuing},
and  
Feller \cite[sec. 6.9]{Feller71v2} note the way CUSUM processes describe queues and storages. 
The amount of a stored commodity increments by the difference between supplied and spent amounts, and it becomes 0 when the demand exceeds
the supply. Similarly, the number of passengers at a bus stop (or the number of jobs in a queue) changes by the difference between 
the arrived and departed ones unless a bus takes all. Similarly, the waiting time of a job in a queue equals the sum of the remaining waiting time 
and the service time of the previous job; that is unless the queue is empty at the time of the job's arrival, in which case its waiting time is zero.
After \cite{Lindley1952queuing}, CUSUM processes have also been called queuing processes or Lindley processes, especially in the queuing theory and probability literature.

CUSUM processes have been broadly used in quality and process control, for monitoring the ``cumulative increase over background''
\cite{ChangGan01, Gibbons99, Lucas82, Woodall85}. 
(Our definition of the CUSUM process should not be confused with plain sums or sums of deviations from the mean that are sometimes called CUSUM
charts \cite{HawkinsOlwell12}. 
To distinguish, \cite{Montgomery09} and \cite{SibandaSibanda07}   call \num{CUSUMdefinition} an upward or upper CUSUM chart.)

Ewan Page 
\cite{Page54} introduced CUSUM processes for sequential change-point detection and showed that repeated testing for change-points
results in a stopping rule that detects a change when the CUSUM process $W_n$ crosses a constant threshold. Later, 
\cite{Lorden71,Moustakides86,Ritov90} proved optimality of Page's CUSUM change-point detection algorithm with respect to various forms 
of optimality. Lorden 
\cite{Lorden71} and Moustakides 
\cite{Moustakides86} showed that it minimizes the worst-case expected delay 
under a fixed rate of false alarms, while Ritov 
\cite{Ritov90} demonstrated the Bayesian nature of the procedure, under a special sequentially
evolving prior distribution which also appears to be the least favorable. See \cite[sec. 2.2]{BassevilleNikiforov93} and 
\cite[sec. 8.2]{TartakovskyNikiforovBasseville14} for an overview. 

In this article, we explore and establish new properties of the CUSUM process, its running maximum, and its moments that find direct applications in change-point detection problems and queuing systems.

We start with the classical approach of Frank Spitzer 
\cite{Spitzer66} that allows studying the marginal and joint distributions of a random walk and its maximum.
Equations (1.3, 6.5) of \cite{Spitzer1956combinatorial} and chapter 17 of \cite{Spitzer66} develop power series of the
corresponding characteristic functions expressed in terms of the moments of rectified partial sums. These rectified sums
\[
 S_n^+ = \max\l S_n, 0\r,  \mbox{ where } S_n = Y_1+\ldots+Y_n,
\]
are simpler variables than $W_n$, and their moments can be derived explicitly or bounded with straightforward expressions. 
Here and later, we use the notation $x^+=\max\l x,0\r$.

In particular, we apply the method of Spitzer to devise recursive computational algorithms for the moments and the moment generating function (MGF) of CUSUM. We derive lower and upper bounds for the CUSUM's MGF, a function that is difficult to study by simulations due to a high variance of the CUSUM's sample MGF.
We show that this MGF has a linear asymptote at a special point $\lam^\ast$ which separates three different types of the
asymptotic behavior of the moment generating function of the CUSUM process. This result is counterintuitive, because the CUSUM process 
renews every time when it hits zero which would suggest its eventual stationarity.  

Direct applications
of these results are in the detection of single and multiple transient change-points as well as the queuing systems. 
Transient changes occur when an observed sequence changes abruptly from an initial distribution $F$ to a new distribution $G$
and then returns to the distribution $F$ after an unknown period of time. As shown in \cite{BaronMalov23transient}, 
$\Lambda_n = \max_{[0,n]} W_t$
is the likelihood-ratio test statistic for the occurrence of a transient change during the interval $[0,n]$, where the CUSUM
process $W_t$ is built upon log-likelihood ratios $Y_i = (dG/dF)(X_i)$ of the observed variables $X_i$. Sequentially,
a transient change is detected at a stopping time when $\Lambda_n$ crosses a pre-selected threshold $h = h(\alpha)$, 
which has to be determined to control the false alarm probability at a desired level $\alpha$. Calculation of $h(\alpha)$
requires the study of a running maximum $\max_{[0,n]} W_t$, which reduces, by the virtue of Doob's maximal inequality
(\cite{Shiryaev95}, Section VII-3;
 \cite{stroock2013mathematics}, Section 7.1.1), to the estimation of exponential moments of $W_n$. 
 
The Spitzer approach, adapted and developed for the study of CUSUM process in this paper, leads to inequalities for the CUSUM exponential moments. 
Consequently, we obtain upper bounds for the threshold $h(\alpha)$ that appear reasonably tight when compared against
precise thresholds obtained as upper $\alpha$-quantiles of the corresponding test statistics. Moreover, our derived 
thresholds are universal for any base distribution $F$ and disturbed distribution $G$. According to \cite{BaronMalov23transient}, 
section 4, the same threshold will also serve the detection of multiple transient changes, controlling the {\sl familywise} false
alarm rate at the same level $\alpha$ without an additional correction for multiple comparisons.

The paper is organized as follows. Adaptation of the Spitzer approach to CUSUM processes is detailed in Section~\ref{sec:moments}, resulting in exact expressions for the moments and exponential moments of the CUSUM process.
In Section~\ref{sec:recursive}, building up on these results and making use of the recursive properties of Bell polynomials
\cite{bell1927partition,bell1934exponential}, we derive recursive and matrix formulas 
as fast computational tools for these CUSUM moments. Applied to change-point analysis, we note a visibly linear behavior 
of exponential moments, a somewhat surprising observation 
which we confirm in Sections~\ref{sec:Doob} and \ref{sec:linbounds} by proving existence of
linear lower and upper bounds for these moments, and ultimately, existence of a linear asymptote. 
Section~\ref{subsec:general} characterizes three types of behavior of the MGF of the CUSUM process, determined by the
mentioned value of $\lamstar$. Applications of the obtained results and thresholds to change-point detection, queuing systems, and biological sequence alignment are discussed in Section~\ref{sec:appl}.  Finally, the Appendix includes a discussion of the application of Spitzer's formula to CUSUM processes. 
To this end, we give an alternative simple proof of the classical Spitzer's formula in \cite{Spitzer1956combinatorial}, 
adapted for CUSUM, demonstrate the failure of Wendell's approach \cite{Wendell58}
to provide a similar framework for the running maximum of CUSUM, and consequently, formulate a conjecture that a similar
expression for the maximum of CUSUM does not exist in a closed form.

\section{Spitzer's approach and exact computation of CUSUM's moments}
\label{sec:moments}

In this section, the classical Spitzer's method of generating functions is extended to CUSUM processes, leading to exact expressions for their moments and the moment generating function.

Let $Y_1,Y_2,\ldots$ be independent and identically distributed (i.i.d.) random variables, $S_0=0$  and $S_n=\sum_{i=1}^n Y_i$ be their partial sums, and $R_n=\max_{[0,n]}S_k$ be the corresponding
sequential maximal values, $n\in\Nb$. The CUSUM process is defined as follows:
\begin{equation}
\label{equ:CUSUM}
W_0=0; \quad W_n=S_n-\min_{[0,n]} S_k.
\end{equation}
This definition is equivalent to the recursive formula \num{CUSUMdefinition}. It is clear that 
\be
W_n=\max_{k\in [0,n]}(S_n-S_k),
\nnn
and therefore, $W_n$ has the same distribution as $R_n$. 
Also, $W_n$ admits a recursive representation \num{CUSUMdefinition}, which is not difficult to derive from \num{equ:CUSUM}.

A classical approach by Frank Spitzer  
\cite{Spitzer66}
  allows to study the distribution of a random walk and its running maximum. 
  We apply a similar methodology to the cumulative-sum process (CUSUM) and obtain
  exact expressions for its moments and exponential moments.
  
  The Spitzer's formula, originally published in \cite{Spitzer1956combinatorial}, yields
\begin{equation}
\label{equ:spz2}
\sum_{n=0}^{\infty} t^n \phi_n(u+v,-v)=\exp\Bigl[ \sum_{n=1}^{\infty}\frac{t^n}{n}(\psi^{+}_n(\lam)+\psi^{-}_n(v)-1)\Bigr],
\end{equation}
and, in the particular case of $v=0$,
\begin{equation}
\label{equ:spz}
\sum_{n=0}^{\infty} t^n \phi_n(u) =\exp\Bigl[ \sum_{n=1}^{\infty}\frac{t^n}{n}\psi^{+}_n(u)\Bigr],
\end{equation}
where $\phi_n(u,v)=\E e^{i(u W_n+v S_{k:n})}$, $\psi^{+}_n(u)=\E e^{iuS_n^+}$, $\psi^{-}_n(v)=\E e^{ivS_n^-}$, $\phi_n(u)=\E e^{iuW_n}$  are characteristic functions of the corresponding random variables
and vectors, $S_{k:n}=S_n-S_k$, and $S_n^+ = \max\l S_n, 0\r$ and $S_n^- = -\min\l S_n, 0\r$ are rectified random walks. 
The MGF can be obtained in a similar manner, as
\begin{equation}
\label{equ:mgf}
\sum_{n=0}^{\infty} t^n\E e^{\lam W_n}=\exp\Bigl[ \sum_{k=1}^{\infty} \frac{t^k}{k}\E e^{\lam S^+_n}\Bigr].
\end{equation}

The most famous application of Spitzer's identity is a convenient representation for the mean value.
We extend this result to the moments and exponential moments of the CUSUM process.

\begin{theorem}[Exact expressions for the moments and MGF]\label{thm:moments}
For the CUSUM process $\{W_n,\, n\ge 0\}$, the running mean, variance, and MGF, if they exist,
admit the following exact expressions:
\be
\E W_n&=& \d \sum_{k=1}^n \frac{1}{k} \E S_k^+, \label{meanW} \\
\D W_n&=& \d \sum_{k=1}^n \frac{1}{k} \D S_k^+ +\sum_{k=1}^n \frac{1}{k} (\E S_k^+)^2-\sum_{{1\le k_1\!,k_2\le n}: \atop k_1+k_2>n} \frac{1}{k_1k_2} \E S_{k_1}^+\E S_{k_2}^+, \label{varianceW}\\
\E e^{\lam W_n}&=& \d \sum_{(k_1,\ldots,k_n):\atop\sum_{j=1}^n jk_j=n} \prod_{r=1}^n \frac{(\E e^{\lam S^+_r})^{k_r}}{r^{k_r}k_r!}. \label{expW}
\eeq
\end{theorem}

\begin{proof}
(i) 
Differentiating (\ref{equ:spz}) with respect to $u$ and setting $u=0$, we obtain
\Be\label{SpitzerDiff}
\sum_{n=0}^{\infty} \E W_n t^n =(1-t)^{-1} \sum_{k=1}^{\infty} \E \Bigl(\frac{t^k}{k} S_k^+ \Bigr)=\sum_{j=0}^{\infty} t^j \sum_{k=1}^{\infty} \E \Bigl(\frac{t^k}{k} S_k^+ \Bigr).
\Eeq
Equating coefficients of the linear terms in the left and the right hand sides of \num{SpitzerDiff}
brings \num{meanW}.
\\[4mm]
(ii) 
Differentiating (\ref{equ:spz}) with respect to $u$ twice and setting $u=0$, we have that 
\be
\sum_{n=0}^{\infty} \E W^2_n t^n =(1-t)^{-1}\Bigl[ \Bigl(\sum_{k=1}^{\infty} \frac{t^k}{k} \E S_k^+ \Bigr)^2+\sum_{k=1}^{\infty} \frac{t^k}{k} \E  (S_k^+)^2 \Bigr],
\nnn
and therefore,
\be
\sum_{n=1}^{\infty} t^n\E (W^2_n-W^2_{n-1})=\Bigl(\sum_{k=1}^{\infty} \frac{t^k}{k} \E S_k^+ \Bigr)^2+\sum_{k=1}^{\infty} \frac{t^k}{k} \E  (S_k^+)^2.
\nnn
Note that 
\be
\Bigl(\sum_{k=1}^{\infty} \frac{t^k}{k} \E S_k^+ \Bigr)^2=\sum_{k_1=1}^{\infty}\sum_{k_2=1}^{\infty} t^{k_1+k_2}m^+_{k_1} m^+_{k_2}
=\sum_{k=2}^{\infty} t^k \Bigl[\sum_{k_1\!,k_2\geq 1: \atop k_1+k_2=k} m^+_{k_1} m^+_{k_2}\Bigr],
\nnn
where $m^+_k=\E S^+_k/k$, $k\in\Nb$. 
Equating the coefficients of different powers of $t$ yields $\E W_1^2=\E (S^+_1)^2=\E (Y^+_1)^2$ and
\be
\E (W^2_n-W^2_{n-1})=\frac{\E  (S_n^+)^2}{n}+\sum_{k_1\!,k_2\geq 1: \atop k_1+k_2=n} m^+_{k_1} m^+_{k_2}, \;\mbox{for}\; n\geq 2.
\nnn
Then, we obtain recursively that 
\be
\E W^2_n=\sum_{k=1}^n \frac{1}{k} \E (S_k^+)^2+\sum_{k_1\!,k_2\geq 1: \atop k_1+k_2\leq n} m^+_{k_1} m^+_{k_2}.
\nnn
Applying \num{meanW} results in the expression,
\be
(\E W_n)^2=\Bigl(\sum_{k=1}^n m^+_k\Bigr)=\sum_{k_1=1}^n\sum_{k_2=1}^n m^+_{k_1} m^+_{k_2}.
\nnn
Hence,
\begin{eqnarray*}
\D W_n&=&\E W^2_n-(\E W_n)^2=\sum_{k=1}^n \frac{1}{k} \E (S_k^+)^2-\!\!\!\!\sum_{{1\le k_1\!,k_2\le n}: \atop k_1+k_2>n} m^+_{k_1} m^+_{k_2}
 \nn &=&
\sum_{k=1}^n \frac{1}{k} \D S_k^+ +\sum_{k=1}^n k (m^+_k)^2-\!\!\!\!\sum_{{1\le k_1\!,k_2\le n}: \atop k_1+k_2>n} m^+_{k_1} m^+_{k_2},
\end{eqnarray*}
and \num{varianceW} is proved.
\\[4mm]
(iii) 
The $d$-th derivative of the left-hand side of \num{equ:mgf} at $t=0$ is equal to $d! \E e^{\lam W_d}$,
which is the key for expressing the exponential moment $\E e^{\lambda W_d}$. Introduce functions $g_{\lam}(t)=\sum_{k=1}^{\infty} \frac{t^k}{k}\E e^{\lam S^+_k}$ for any $\lam$ from the domain of MGF of $S_k^+$ for all $k\ge 1$, which coincides with the domain of MGF of $Y_1^+$. Then, 
differentiate the right-hand side of \num{equ:mgf} according to the Fa\`{a} di Bruno's formula (e.g., (2.5) in \cite{lukacs1955applications}), we obtain the $d$-th derivative
\begin{eqnarray*}
\exp(g_\lam (t))^{(d)}&=&\exp(g_\lam (t))\!\!\sum_{(m_1,\ldots,m_d):\atop\sum_{j=1}^d jm_j=d}\!\frac{d!}{m_1!\cdots m_d!} \prod_{r=1}^d\Bigl(\frac{g_\lam ^{(r)}(t)}{r!}\Bigr)^{m_r}.
\end{eqnarray*}
Taking into account that $g_\lam (0)=0$ and $g_\lam ^{(r)}(0)=(r-1)!\, \E e^{\lam S^+_r}$, we have  
\begin{eqnarray}
\E e^{\lam W_d}&=& \left. \frac{1}{d!} \exp(g_\lam (t))^{(d)} \right|_{t=0} 
 = \sum_{(m_1,\ldots,m_d):\atop\sum_{j=1}^d jm_j=d}\prod_{r=1}^d\frac{(\E e^{\lam S^+_r})^{m_r}}{r^{m_r}m_r!},\nonumber
\end{eqnarray}
which completes the proof of \num{expW}.

\end{proof}

Note that formulas (\ref{meanW}) and (\ref{varianceW}) were obtained in \cite{AnisGharib80} in more general case of exchangeable random variables without using Spitzer's identity. Taking further derivatives of \num{equ:spz}, one can obtain exact expressions for higher moments
of the CUSUM process. The exponential moment  (if exists) can be obtained from the exponential moments of $S^+_n$ by 
substituting $\lam=1$ in (\ref{expW}):
\begin{equation}
\E e^{W_n}=\d \sum_{(k_1,\ldots,k_n):\atop\sum_{j=1}^n jk_j=n} \prod_{r=1}^n \frac{(\E e^{S^+_r})^{k_r}}{r^{k_r}k_r!}.\nonumber
\end{equation}

An alternative proof of Spitzer's identities (\ref{equ:spz}) 
and  (\ref{equ:mgf}) is given in Appendix  \ref{appA}. 
We offer a conjecture that Spitzer's identity for 
the maximum of a CUSUM process does not exist in a closed form. 
Implicit arguments in support of this conjecture are in Appendix  \ref{appB}.

\section{Bell's polynomials and recursive computation of CUSUM's 
 moments}\label{sec:recursive}

In this section, we use the results of Theorem~\ref{thm:moments} and recursive properties of Bell's polynomials
to derive fast computational tools for the moments and MGF of CUSUM processes. Recursive and matrix
methods are proposed for this purpose.

Theorem~\ref{thm:moments} expresses  moments of the CUSUM process $W_n$ in terms of a simpler stochastic
process $S_n^+$. However, the right-hand sides of (\ref{meanW}--\ref{expW}) are still cumbersome 
for efficient computation. In this section, we derive a recursive algorithm for the fast and efficient computation of the first moments and MGF of CUSUM. Recursive expressions exist for higher moments of the CUSUM too, and they can be derived similarly or directly from Theorem~\ref{thm:moments}.

\begin{theorem}[Recursive computation of CUSUM's moments and MGF]\label{thm:recursive}
The expected value $\M_n=\E W_n$, the variance $\V_n=\D W_n$, and the MGF $M_n(\lam) = \E e^{\lam W_n}$ of a CUSUM process at time $n$ satisfy the following recursive relations,
\begin{eqnarray}
&&\l\ba{lll}
\M_0 &=& 0, \\
\M_{n} &=& \d \M_{n-1} + \frac{1}{n} \E S_n^+ \mbox{, for } n\ge 1;
\ea\right. \label{meanWrecursiveThm}
\\[1mm]
&&\l\ba{lll}
\V_0 &=& 0, \\
\V_{n} &=& \d \V_{n-1} + \frac{1}{n} \D  S_n^+   + \sum_{k=1}^n \frac{1}{k(n+1-k)} \E  S_k^+  \E  S_{n+1-k}^+ \mbox{, for } n\ge 1;
\ea\right. \label{varianceWrecursiveThm}
\\[1mm]
&& \l\ba{lll}
M_0(\lam ) &=& 1 \\
M_{n+1}(\lam ) &=& \d\frac{1}{n+1} \sum_{k=0}^n M_k(\lam ) \E e^{\lam  S^+_{n-k+1}} \mbox{, for } n\ge 0  
\ea\right. \label{exponentialWrecursiveThm}
\end{eqnarray}
and any $\lam $ from the domain of the MGF of $Y_1^+$.
\end{theorem}

\begin{proof}
Recursive equations \num{meanWrecursiveThm} and \num{varianceWrecursiveThm} are deduced directly from
\num{meanW} and \num{varianceW}, with obvious arithmetic. The initial condition follows from $W_0=0$. 
Recursive expressions for the CUSUM's MGF are proven below.

Application of the Fa\`{a} di Bruno's formula for the higher order derivatives of a composite function
allows to express the exponential moment of CUSUM in terms of Bell polynomials, 
\cite{bell1927partition,bell1934exponential}. A {\sl complete exponential Bell polynomial} is defined as
\be
B_d(x_1,\ldots,x_d)=\sum_{k=1}^d B_{dk}(x_1,\ldots,x_{d-k+1}),
\nnn
where
\be
B_{dk}(x_1,\ldots,x_{d-k+1})=d! \sum_{(m_1,\ldots,m_d)\in {\mathcal D_{dk}}}\prod_{i=1}^{d-k+1} \frac{x_i^{m_i}}{(i!)^{m_i}m_i!}
\nnn
are
{\sl incomplete exponential Bell polynomials} and ${\mathcal D_{nt}}\!=\!\{(k_1,\ldots,k_n)\!:\!\sum_{j=1}^{n-t+1} \!k_j\!=\!t,\sum_{j=1}^{n-t+1} j k_j\!=\!n\}$.


In terms of Bell polynomials, the CUSUM MGF is expressed as
\begin{eqnarray}
\label{equ:bellexp}
M_d(\lam) &=& \E e^{\lam W_d}=\sum_{(m_1,\ldots,m_d):\atop\sum_{j=1}^d jm_j=d}\prod_{r=1}^d\frac{(\E e^{\lam S^+_r})^{m_r}}{r^{m_r}m_r!}\!=\! \frac{1}{d!}B_{d}\bigl(0!\,\E e^{\lam S^+_1} ,\ldots,(d-1)!\,\E e^{\lam S^+_{d}}\bigr).
\nonumber
\end{eqnarray}

Bell's polynomials admit a recursive formula (e.g., Theorem 1 of \cite{bulo2012efficient})
\be
B_{n+1}(y_1,\ldots,y_{n+1})=\sum_{k=0}^{n} \binom{n}{k} B_{n-k}(y_1,\ldots,y_{n-k})y_{k+1}, 
\nnn
which is not practical for the calculation of $\E e^{\lam W_n}$ for large $n$ because of a very large magnitude of
$\binom{n}{k}$. Instead, we introduce {\sl rescaled Bell polynomials}
\Be\label{rescaledBell}
\widetilde B_n(x_1,\ldots,x_n)=\frac{1}{n!} B_n\bigl(0!\,x_1,\ldots,(n-1)!\,x_n\bigr).
\Eeq
They are directly connected to the CUSUM's MGF,
\Be\label{EexpW_Bell}
\E e^{\lam W_n}=\widetilde B_n\bigl(\E e^{\lam S^+_1},\ldots,\E e^{\lam S^+_n}\bigr),
\Eeq
for all $d\in\Nb$. The following recursive relation allows fast computation of $\E e^{\lam W_n}$ from known MGFs
$\E e^{\lam S^+_k}$, $k\in\Nb$, even when $n$ is large:
\begin{eqnarray}
\label{equ:BellRec}
\lefteqn{\widetilde B_{n+1}(x_1,\ldots,x_{n+1}) = \frac{1}{(n\!+\!1)!}\sum_{k=0}^n C_n^k B_{n-k}\bigl(0!\,x_1,\ldots,(n\!-\!k\!-\!1)!\,x_{n-k}\bigr)k!\,x_{k+1}} \\
\ \ &=& \frac{1}{(n+1)!}\sum_{k=0}^n \frac{n!k!(n-k)!}{k!(n-k)!} \widetilde B_{n-k}\bigl(x_1,\ldots,x_{n-k}\bigr)\,x_{k+1}\nonumber\\
\ \ &=&\frac{1}{(n+1)}\sum_{k=0}^n \widetilde B_{n-k}\bigl(x_1,\ldots,x_{n-k}\bigr)\,x_{k+1}
= \frac{1}{(n+1)}\sum_{k=0}^n \widetilde B_{k}\bigl(x_1,\ldots,x_{k}\bigr)\,x_{n-k+1}. \nonumber
\end{eqnarray}

The recursive step is then obtained from \num{EexpW_Bell} and \num{equ:BellRec},
%
\be
\lefteqn{M_{n+1}(\lam ) = \widetilde B_{n+1}\bigl(\E e^{\lam S^+_1},\ldots,\E e^{\lam S^+_{n+1}}\bigr)}
 \nn  && = \frac{1}{(n+1)}\sum_{k=0}^n \widetilde B_{k}\bigl(\E e^{\lam S^+_1},\ldots,\E e^{\lam S^+_k}\bigr)\,\E e^{\lam S^+_{n-k+1}} 
        = \frac{1}{(n+1)}\sum_{k=0}^n M_k(\lam ) e^{\lam S^+_{n-k+1}}.
\nnn
\end{proof}

Implementation of the recursive algorthm provided in Theorem~\ref{thm:recursive} requires computation of
exponential moments of rectified sums $x_k^{(\lam)} = \E e^{\lam S_k^+}$. An example of such computation is given in 
\num{ESplus}. After that, the CUSUM's 
exponential moments $M_n{(\lam)}=\E e^{\lam W_n}$ for $n=0,1,2,\ldots,$ can be calculated as a solution of a matrix equation or 
as a result of iterations.

\subsection{Matrix computation of exponential moments}

From Theorem~\ref{thm:recursive}, denoting the entire vector of 
\[
\l M_{k}(\lam) = M_{k}(\lam;x_1^{(\lam)}\ldots,x_k^{(\lam)}), \, k=0,\ldots,n\r
\]
by ${M(\lam)}$, we have the following matrix equation for ${M(\lam)}$,
\[
M(\lam) = A(\lam)M(\lam) + e,
\]
where 
%
{\arraycolsep=3.5pt
\Be\label{comp_matrix}
M(\lam) = \left(\ba{c}M_0(\lam)\\ M_1(\lam)\\ M_2(\lam)\\ \cdots \\ M_n(\lam) \ea\right)\!, \
A(\lam) = \left(\ba{ccccc}
  0\!&\!0\!&\!0\!&\!\cdots\!&\!0 \\
  x_1^{(\lam)}\!&\!0\!&\!0\!&\!\cdots\!&\!0 \\[1mm]
  \frac{1}{2}x_2^{(\lam)}\!&\!\frac{1}{2}x_1^{(\lam)} \!&\!0\!&\!\cdots\!&\!0 \\[2mm]
  \frac{1}{3}x_3^{(\lam)}\!&\!\frac{1}{3}x_2^{(\lam)}\!&\!\frac{1}{3}x_1^{(\lam)} \!&\! \cdots\!&\!0 \\ 
  \cdots\!&\!\cdots\!&\!\cdots\!&\!\cdots\!&\!\cdots \\
  \frac{1}{n}x_n^{(\lam)}\!&\!\frac{1}{n}x_{n-1}^{(\lam)}\!&\! \frac{1}{n}x_{n-2}^{(\lam)} \!&\! \cdots\!&\!0 \\  
  \ea\right)\!, \
 e = \left(\ba{c}1\\ 0\\ 0\\ \cdots \\ 0 \ea\right)\!.
\Eeq
}
 Thus, we get an exact matrix solution for all exponential moments of the CUSUM process up to time $n$, 
%
 \[
 M(\lam) = \left(\!\!
 \ba{c} \E e^{\lam W_0} \\ \E e^{\lam W_1} \\ \E e^{\lam W_2} \\ \cdots \\ \E e^{\lam W_n} \ea \!\!\right)
  =
 (I-A(\lam))^{-1}e = \left(\!\!\ba{ccccc}
  1\!\!&\!\!0\!\!&\!\!0\!\!&\!\!\cdots\!\!&\!\!0 \\
  -x_1^{(\lam)}\!\!&\!\!1\!\!&\!\!0\!\!&\!\!\cdots\!\!&\!\!0 \\[1mm]
  -\frac{1}{2}x_2^{(\lam)}\!\!&\!\!-\frac{1}{2}x_1^{(\lam)} \!\!&\!\!1\!\!&\!\!\cdots\!\!&\!\!0 \\[2mm]
  -\frac{1}{3}x_3^{(\lam)}\!\!&\!\!-\frac{1}{3}x_2^{(\lam)}\!\!&\!\!-\frac{1}{3}x_1^{(\lam)} \!\!&\!\! \cdots\!\!&\!\!0 \\ 
  \cdots\!\!&\!\!\cdots\!\!&\!\!\cdots\!\!&\!\!\cdots\!\!&\!\!\cdots \\
  -\frac{1}{n}x_n^{(\lam)}\!\!&\!\!-\frac{1}{n}x_{n-1}^{(\lam)}\!\!&\!\! -\frac{1}{n}x_{n-2}^{(\lam)} \!\!&\!\! \cdots\!\!&\!\!1 \\  
  \ea\!\!\right)^{-1}
  \!\!\!\!\!\!
  \left(\!\!\ba{c}1\\ 0\\ 0\\ \cdots \\ 0 \ea\!\!\right).
 \]
The corresponding algorithms appear rather fast and computationally cost-efficient, as we observe in the following example.

\subsection{Recursive computation of exponential moments}\label{subsec:expmoments}

Theorem \ref{thm:recursive} also suggests an iterative computation of CUSUM moments and MGF.
The corresponding algorithms appear rather fast and computationally cost-efficient, as we observe in the following example.


\begin{figure}
\begin{picture}(0,140)(0,0)
\put(-130,-22){\scalebox{1}[0.9]{\includegraphics[scale=.54]{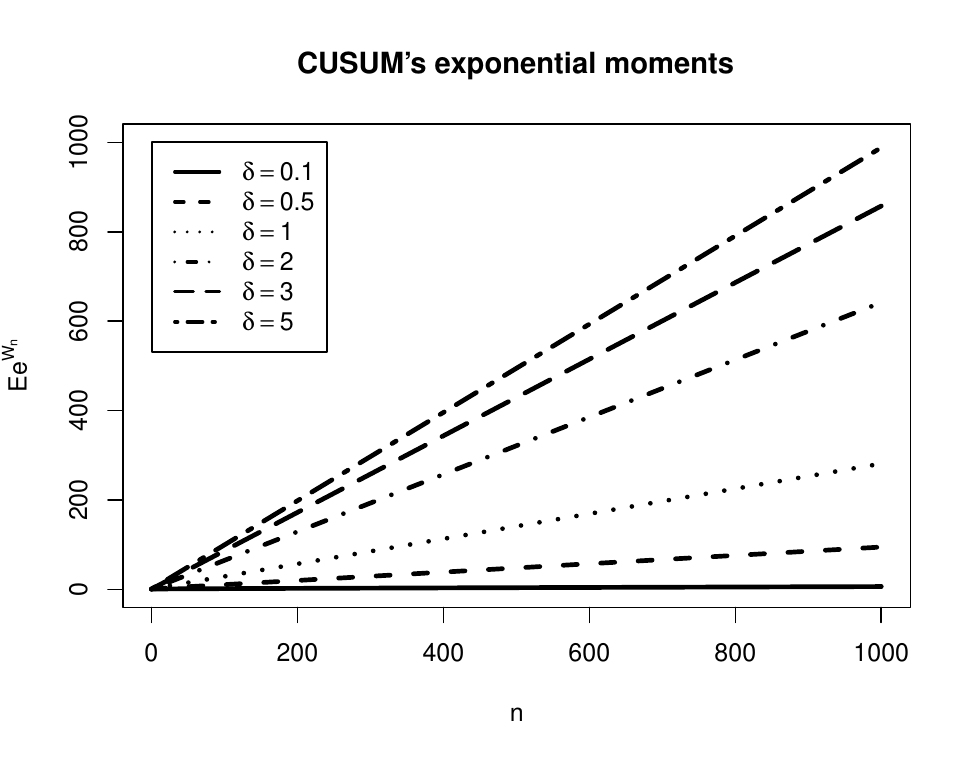}}}
\end{picture}
\caption{CUSUM's exact exponential moments for various detectable differences $\delta$.\label{fig:expmoments}}
\end{figure}

A CUSUM process is often used for the 
fast change-point detection, representing the optimal non-Bayesian detection procedure,
according to several standard criteria 
\cite{borovkov1999asymptotically,LaiIEEE98,Lorden71,Moustakides86,Page54}, 
Chap. 2.2, 5.2 of \cite{BassevilleNikiforov93},
Chap. 8.2 of \cite{TartakovskyNikiforovBasseville14}. 
Properties of this CUSUM process
provide the information needed for the computation of a detection threshold that satisfies the
desired control of the rate of false alarms.

For the detection of an abrupt change from a {\sl default} distribution $F$ with density $f$ to a {\sl disturbed} distribution $G$ with density $g$, with respect to an underlying measure $\mu$, the CUSUM process is generated by log-likelihood ratios
\[
Y_t = \log\frac{g(X_{t+1})}{f(X_{t+1})}
\]
and the associated random walk $S_n$. Exponential moments of rectified sums in this case are
expressed in terms of the original random walk as
\be
x_n^{(1)} &=& \E_F e^{S^+_n} = \E_F \Bigl(e^{S_n}\indi_{\{S_n\geq 0\}}+\indi_{\{S_n < 0\}}\Bigr)
\nn &=& \E_F \prod_{i=1}^n\frac{g(X_i)}{ f(X_i)}\indi_{\{S_n\geq 0\}}+\P(S_n<0) = {\Pr}_G(S_n\geq 0)+\P_F(S_n<0). \label{ESplus}
\eeq
Two probabilities in \num{ESplus} are calculated directly in many statistical software languages,
when the distributions $F$ and $G$ belong to a convolution-closed family, such as Normal, Gamma, Chi-Square, Poisson, Binomial, or Negative Binomial. 

For example, for the detection of a change from the Normal$(\theta_0,\sigma)$ distribution $F$ with
mean $\theta_0$ and standard deviation $\sigma$ to the Normal$(\theta_1,\sigma)$ distribution $G$,
the partial sums $S_n$ are Normal$(-\delta^2/2, \delta)$, where $\delta=|\theta_1-\theta_0|/\sigma$ 
is the standardized detectable difference, which also represents {\sl the magnitude of change}. 

Algorithm \num{exponentialWrecursiveThm} for the computation of the CUSUM's first $N$ exponential moments is 
then coded in a few lines in R software.

\vskip -3mm

\bc{\scalebox{0.88}{            
 \tt{\bt{ll}
 M <- numeric(N+1);                        & \# First element is M[0], then M[1], etc. \\
 M[1] <- 1;                                & \# EU[1] = E$e^{W[0]}$ = 1           \\
 x <- 2*pnorm(delta*sqrt(1:N)/2);          & \# Exponential moments of rectified sums  \\
 for (n in 1:N)\{                          & \# Iterative computation of CUSUM's    \\
 \ \ \   M[n+1] <- sum(M[n:1] * x[1:n])/n; & \# \ \ \ \ \ \ \ \ \ exponential moments                 \\
 \ \ \     \}
\et}}}\ec

\vskip 2mm


For different magnitudes $\delta$ from a standardized difference of 0.1, which is difficult to detect, to an easily
detectable difference of 5.0 standard deviations, the CUSUM's calculated exponential moments are depicted in
Figure~\ref{fig:expmoments}. Noticeably, their behavior is linear for any $\delta$.

\begin{figure}
\begin{picture}(0,130)(240,0)
\put(35,-18){\scalebox{.31}[.33]{\includegraphics[scale=.95]{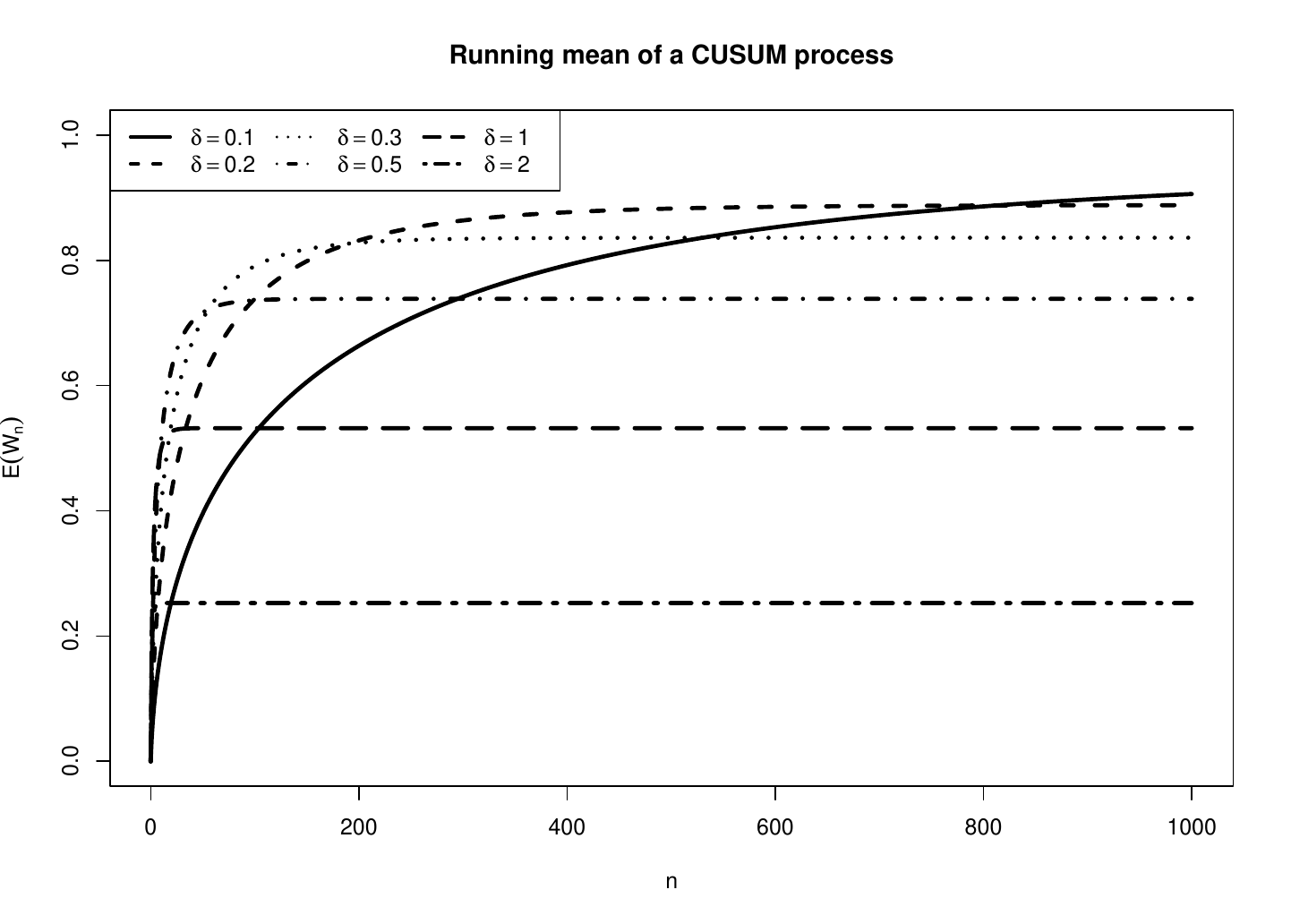} }}
\put(240,-18){\scalebox{.31}[.33]{\includegraphics[scale=.95]{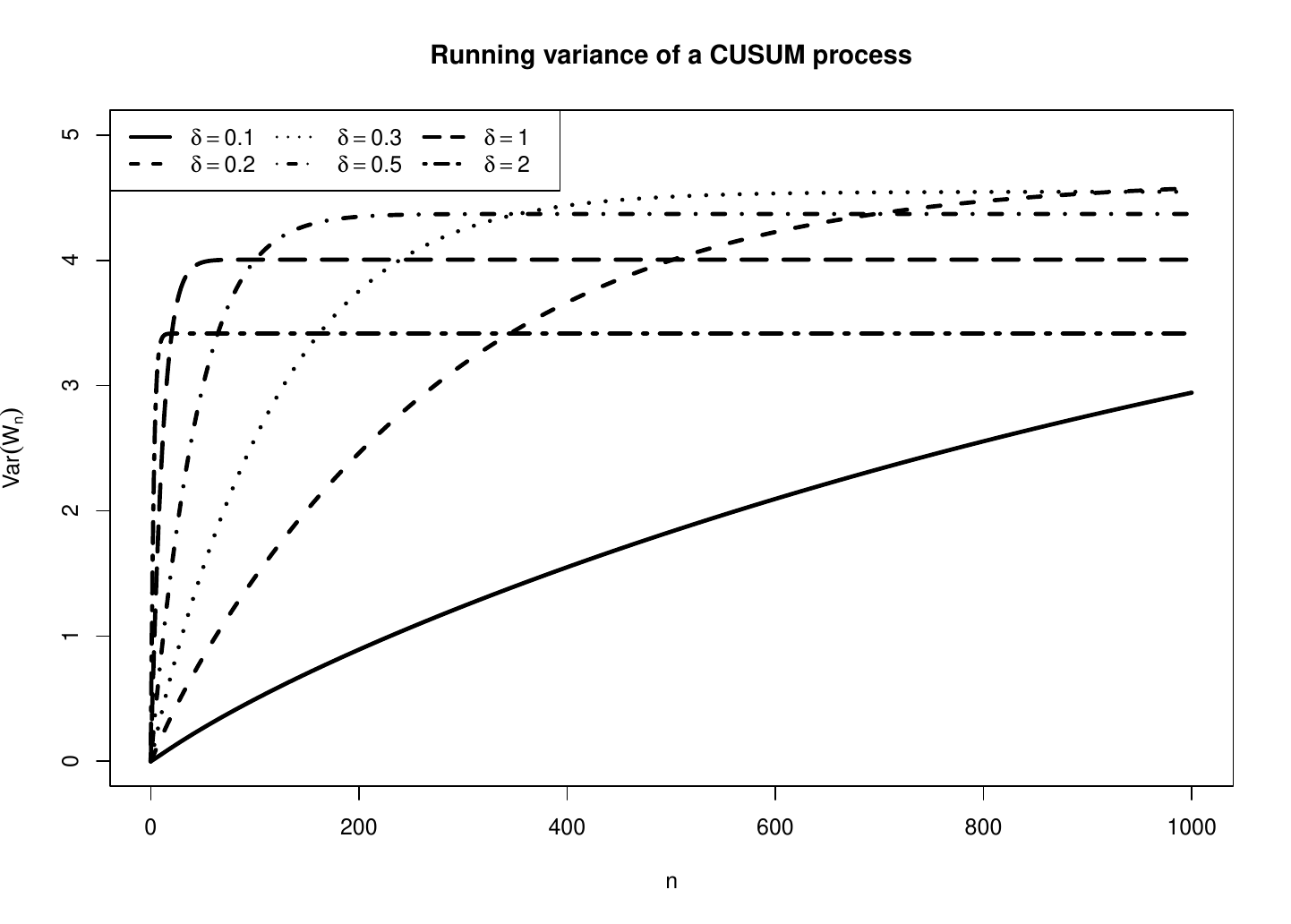}  }  }
\end{picture}
\caption{CUSUM's mean and variance for various detectable differences $\delta$. \label{fig:mean_and_variance}}
\end{figure}

This is somewhat counterintuitive. We know that the CUSUM process $W_n$ renews at every moment $n$ when $W_n=0$, 
which happens infinitely often due to a negative drift $(-\delta^2/2)$ of random walk $S_n$, at i.i.d. intervals. Then, one would characterize process $W_n$ as a sequence of i.i.d. segments. Thus, it seems reasonable to expect the distribution of $e^{W_n}$ to stabilize.
Indeed, this is the case with the expected value $\E W_n$ and variance $\D W_n$, as seen on Figure~\ref{fig:mean_and_variance}, where these moments
are calculated according to the recursive formulas in Theorem~\ref{thm:recursive}.
Yet, we observe in Figure~\ref{fig:expmoments} that instead of tending to a constant, the expected value
$\E e^{W_n}$ has a precise linear trend. This phenomenon is explored and proven in the next two Sections.

Further discussion and application of our probability results to change-point analysis are detailed in Section~\ref{sec:appl}.

\section{Maximal inequalities for CUSUM}
\label{sec:Doob}

Sections~\ref{sec:Doob} and \ref{sec:linbounds} provide the theoretical base for the asymptotic behavior of exponential
moments and MGF of CUSUM, $M_n(\lam)=\E e^{\lam W_n}$, including their linearity for certain values of $\lam$. 
First, we derive a reasonably sharp linear upper bound for these moments from the Doob’s maximal inequality 
and the Doob–Meyer decomposition of submartingales. A linear lower bound for the same moments is then obtained
from the study of segments of random walk $S_n$ and a subsequent use of the Cramer’s large deviations theorem.
Existence of a slant asymptote of exponential moments $\E e^{W_n}$ is proven in Section~\ref{sec:linbounds}.

%

For convenience, introduce a random variable $Y$ with the same distribution as variables $Y_1,Y_2,\ldots$ 
that generate the random walk $S_n$ and the corresponding CUSUM process $W_n$.
At first, we assume existence of the exponential moment 

\vskip -4mm

 \Be\label{mgfcondition}
 \E e^{Y}=1.
 \Eeq
In Section~\ref{subsec:general}, we relax this condition replacing it with 
existence of a positive root of the equation $\E  e^{\lambda Y}=1$, a key quantity that identifies and separates
three types of the asymptotic behavior of the CUSUM's exponential moments. 
We note that either condition implies a positive root of the equation
$m(\lambda) = 1$, where $m(\lambda) = \E e^{\lambda Y}$ is the moment generating function of $Y$. Since $m(0)=1$ is the other root, $m'(0)=\E Y$, and $m(\lambda)$ is convex, we conclude that 
\Be\label{negEY}
\E Y < 0.
\Eeq
Preparing for the use of Doob's maximal inequality, consider the process $U_n = e^{ W_n }$, which admits 
the recursive equation 
\Be\label{Urecurs}
U_0 = 1, \ U_{n+1} = \max\{ 1, \, U_n e^{Y_{n+1}}\}
\Eeq
that follows from \num{CUSUMdefinition}. Taking expected values, we obtain 
\[
\E \{ U_{n+1} \,|\, U_1,\ldots,U_n \} \ge U_n \E  e^{Y_{n+1}} = U_n
\]
from \num{mgfcondition}, showing that $U_n$ is a {\sl submartingale} (unlike the CUSUM $W_n$ itself) with respect to the natural filtration $\Fb=\{\Ff_n,n=0,1,\ldots\}$. 

Then, the Doob's maximal inequality can be applied to the process $\{
U_n\}$, resulting in the following inequality,
\be
\Pr \{ \max_{[0,n]} W_t \ge h \}
 = \Pr \{ \max_{[0,n]} U_t \ge e^h \}\le e^{-h}\E (U_n)
 = e^{-h}\E (e^{W_n}).
\label{equ:DoobCUSUM}
 \eeq

\subsection{Upper bounds for exponential moments and tail probabilities}

Next, introduce the sequence $\{\widetilde U_n\}_{n=1}^{\infty}$ defined as the {\sl compensator} of $U_n$ 
(see the Doob-Meyer decomposition, \cite[Chap. VII]{Doob1953}).
Being a compensator means that $\widetilde U_n$ is such a process that $(U_n-\widetilde U_n)$ is a martingale, and it is determined explicitly as
 \be
\widetilde U_1 &=& \E e^{Y^+_1}, \nn
\widetilde U_n &=& \E e^{Y^+_n}+\sum\nolimits_{k=2}^n \E (U_k-U_{k-1}|\Ff_{k-1}) 
 \nn
&=& \E e^{Y^+_1}\!+\!\sum\nolimits_{k=2}^n U_{k-1}(\E  (e^{Y_k\vee (-W_{k-1})}|\Ff_{k-1})-1).
\nnn
By the independence property of $Y_n$ and $\Ff_{n-1}$ for each $n\in\Nb$,
\begin{eqnarray*}
\lefteqn{\E  \bigl(e^{Y_n\vee (-W_{n-1})}|\Ff_{n-1}\bigr) = \E  \bigl(e^{Y_n\!\vee\! (-W_{n-1})}|W_{n-1}\bigr)
}
 \nn
&=& \E \bigl(e^{Y_n}\indi_{\{Y_n>-W_{n-1}\}}|W_{n-1}\bigr)
+e^{-W_{n-1}}\Pr\bigl(Y_n\leq -W_{n-1}|W_{n-1}\bigr) 
 \nn
&=& 1+\E \bigl((e^{-W_{n-1}}-e^{Y_n})\indi_{\{Y_n\leq -W_{n-1}\}}|W_{n-1}\bigr). 
\end{eqnarray*}
Hence,
\begin{eqnarray}
\widetilde U_n&=&\widetilde U_1+\sum\nolimits_{k=2}^n e^{W_{k-1}} \E \bigl((e^{-W_{k-1}}-e^{Y_k})\indi_{\{Y_k\leq -W_{k-1}\}}|W_{k-1}\bigr)\nonumber\\
\label{equ:compensator}
&=&\E  e^{Y_1} +\sum\nolimits_{k=1}^n  \E  \bigl((1-uZ)\indi_{\{Z\leq 1/u\}}\bigr)\bigl|_{u=U_{k-1}},
\end{eqnarray}
where $Z=e^{Y}$, $W_0=0$, and $U_0=1$. Note that $Q(u)= \E \bigl((1-u Z)\indi_{\{Z\leq 1/u\}}\bigr)$ is 
a non-increasing function of $u\in (0,\infty)$, because
\be
Q(v)-Q(u)=-\E  \indi_{\{Z\in (1/v,1/u]\}}+u\E  Z\indi_{\{Z\in (1/v,1/u]\}}-(v-u)\E  Z\indi_{\{Z\leq 1/v\}}\leq 0,
\nnn
for any $0<u<v<\infty$. Therefore, 
\be
Q(U_n)\leq Q(1) = {\E  (1-Z)\indi_{\{Z\leq 1\}} = \E  \left(1-e^{Y}\right)^+}
\nnn
for all $n\in\Nb$. Since
$
\E  U_n=\E  \widetilde U_n
$
by the martingale property of $(U_n-\widetilde U_n)$, we obtain  {\sl the upper bound}
\begin{equation}
\label{equ:marbound}
\E  U_n\leq {1+n\E  \left(1-e^{Y}\right)^+},
\end{equation}
for all $n\in\Nb$. 

{Combining this inequality with \num{equ:DoobCUSUM} yields an upper bound for the 
probability of the CUSUM process to attain threshold $h$ at least once during the time segment $[0,n]$.
\begin{theorem}\label{thm:maxCUSUM_UBgeneral}
For any sample size $n$ and threshold $h$, the tail probability of the current maximum value of the CUSUM process
satisfies the following upper bound, 
\Be\label{maxCUSUM_UBgeneral}
\Pr \{ \max_{[0,n]} W_t \ge h \} \le e^{-h}\left( 1 + n\E (1-e^{Y})^+ \right) \wedge 1.
\Eeq
The CUSUM's exponential moment admits the upper bound
\[
\E e^{W_n} \le {1+n\E  \left(1-e^{Y}\right)^+}.
\]
\end{theorem}
}

\subsection{Lower bounds for exponential moments and tail probabilities}
\label{ssec:lower}

Since $\max_{[0,n]}W_n\geq \max_{[0,n]}Y_k$, the {\sl lower bound} of the CUSUM's MGF can be obtained directly from  (\ref{equ:DoobCUSUM}):
\begin{equation}
\label{equ:GenDoob}
\E  e^{\lambda W_n}\geq \sup_{h>0}\bigl\{ e^{\lambda h}\Pr (\max_{[0,n]} W_t \ge h)\bigr\}
\geq \sup_{h > 0}\bigl\{e^{\lambda h}\Pr (\max_{[0,n]}Y_k \ge h)\bigr\},
\end{equation}
for each $n\in\Nb$.

The lower bound for the tail probability $\Pr(\max_{[0,n]} W_t \ge h)$, and consequently,
for the MGF of $W_n$ can be improved by the following argument. Choose an integer $k$, split the interval $[0,n]$ into
$[n/k]$ subintervals, each of length $k$, and consider the probability for at least one of these segments of 
random walk $S_k = Z_1+\ldots+Z_k$ to exceed threshold $h$. This approach results in the inequality,
\begin{equation}
\label{equ:segmax}
\Pr (\max_{[0,n]} W_t \ge h)\geq 1-\bigl(1-\Pr(S_k \ge h)\bigr)^{[n/k]}.
\end{equation}
Choose $h=ku$ with some positive constant  $u>0$ and let $k=k(n)$ be an {\sl arbitrary} function of $n$, 
taking integer values and increasing to $+\infty$ 
at a rate slower than $n$. Then, by the Cramer's large deviations theorem (e.g., \cite{Dembo09}, Theorem 2.2.3),
the probability in the right-hand side of (\ref{equ:segmax}) equals
%
\be
\Pr (Y_1+\ldots+Y_k \ge ku ) &=& \exp\Bigl(-k\bigl(\inf\{I(x):x\ge u\}+o(1)\bigr)\Bigr)
\nn 
&=& \exp\Bigl(-k\bigl(I(u)+o(1)\bigr)\Bigr),
\label{probdeviation}
\eeq
where $I(x)=\sup_{\lambda\in\Rb} \{\lambda x-\log \E_{F} e^{\lambda Y}\}$ is the Cramer's transformation of the distribution $F$, also known as the {\sl Legendre-Fenchel transformation} of MGF of the distribution $F$, where $o(1)$ tends to $0$ as $n\to\infty$, 
and $\E_{F} e^{\lambda Y}$ can be infinite for some $\lambda$ \cite[Definition 2.2.2]{Dembo09}.

In the application of CUSUM processes to change-point problems, discussed in Sections~\ref{subsec:expmoments} and
\ref{subsec:chp}--\ref{subsec:comparison}, condition \num{mgfcondition} always holds under the original, 
pre-change distribution. Here, we consider two other general cases - the normal distribution and the 
Bernoulli distribution of the original data.

For the normal distribution, $Y\sim N(a,\sigma^2)$, condition \num{mgfcondition} holds under $a=-\sigma^2/2$. 
The Legendre-Fenchel transform for the normal family is $I(x)=(x-a)^2/(2\sigma^2)$. Application of \num{probdeviation} to this case reduces to
\begin{eqnarray}   
\Pr (Y_1+\ldots+Y_k \ge ku ) \nonumber
&=& 1-\Phi\Bigl(\frac{\sqrt k (u-a)}{\sigma}\Bigr)\\
&=&\Phi\bigl(-\frac{\sqrt k (u-a)}{\sigma}\bigr)
\sim\frac{\sigma}{\sqrt{2\pi k}\,|u\!-\!a|} e^{-k\frac{(u-a)^2}{2\sigma^2}},\nonumber
\end{eqnarray}
where 
\be
\frac{1}{\sqrt{2\pi k}\,|u-a|/\sigma}=\exp\Bigl(-k \frac{\frac{1}{2}\log(2\pi k)+\log(|u-a|/\sigma)}{k}\Bigr)=\exp\bigl(-k\, o(1)\bigr).
\nnn
Then $I(x)-x=x^2-2(a+\sigma^2)x+a^2=(x+a)^2\geq 0$, and, therefore, $I(x)\geq x$ for all $x\in\Rb$ with the {\sl unique} root $x_0=-a=\sigma^2/2>0$ of the equation $I(x)=x$.

Under the (transformed) Bernoulli distribution of $Y$ with $\Pr(Y=1)=1-\Pr(Y=-1)=p$, the supremum in the definition of
$I(x)$ is attained at 
%
\[ \lambda=\lambda(x)=\frac{1}{2}\log\frac{(1-p)(1+x)}{p(1-x)},\] and
the Cramer's transformation has the following form,
\be
I(x)=\frac{1-x}{2}\log\Bigl(\frac{1-x}{1-p}\Bigr)+\frac{1+x}{2}\log\Bigl(\frac{1+x}{p}\Bigr)-\log 2.
\nnn
Condition \num{mgfcondition} for the Bernoulli family holds under $p=1/(e+1)$. In this case, the equation $I(x)=x$ again has 
a {\sl unique} solution. Indeed, 
\be
I(x)-x=\frac{1-x}{2}\log\Bigl(\frac{(1-x)(e+1)}{e}\Bigr)+\frac{1+x}{2}\log\bigl((1+x)(e+1)\bigr)-\log 2-x.
\nnn
Differentiating by $x$ yields,
\begin{eqnarray*}
(I(x)-x)'&=&-\frac{1}{2}\Bigl(\log\Bigl(\frac{(1-x)(e+1)}{e}\Bigr)+1\Bigr)+\frac{1}{2}\Bigl(\log\bigl((1+x)(e+1)\bigr)+1\Bigr)-1\\
&=&-\frac{1}{2}\log (1-x)+\frac{1}{2}\log (1+x)-\frac{1}{2}=\frac{1}{2}\log\Bigl(\frac{1+x}{1-x}\Bigr)-\frac{1}{2}.
\end{eqnarray*}
Equating this expression to zero and solving for $x$, we obtain the minimizer of the function $I(x)-x$ at 
\be
x_0=\frac{e-1}{e+1}.
\nnn 
Finally, 
\be
I(x_0)=\frac{1}{1+e}\log(2/e)+\frac{e}{1+e}\log(2e)-\log 2=\frac{e-1}{e+1}=x_0,
\nnn
and, therefore, $x_0$ is a {\sl unique root} of the equation $I(x)=x$. 

The next lemma extends this important property to the general case. 

\begin{lemma}
\label{lem:standard}
Let $\Pr (Y>0)>0$ and $\E  e^{Y}=1$. Then $\inf_{x} (I(x)-x)=0$. Moreover, if $\E  e^{\lambda Y}<\infty$ for some $\lambda>1$,
the equation $I(x)=x$  has a unique root $x_0=\arg\min_{x} (I(x)-x)>0$ that satisfies $\lambda(x_0)=1$, where
$\lambda(x) = \arg\max_{\lambda\in\Rb} \{\lambda x-\log m(\lambda)\}$
is the point of maximum that defines the Legendre-Fenchel transformation of $m(\lambda)$, the MGF of $Y$. 
\end{lemma}

\begin{proof}  
Under the condition $\E  e^{Y}=1$, the moment generating function $m(\lambda)=\E e^{\lambda Y}$ exists for $\lambda\in[0,1]$. Since the Legendre transformation of $\log m(\lambda)$ is a convex function on its domain, 
the function $I^*(x)=I(x)-x=l(x,\lambda(x))-x$ is convex as well. Being the point of maximum of 
$l(x,\lambda)=\lambda x-\log \E  e^{\lambda Y}$ under a fixed $x$, the function $\lambda(x)$ can be obtained from the equation 
\begin{equation}
\label{equ:Legendre}
\frac{\partial}{\partial\lambda} l(x,\lambda)=x-\frac{\E  Ye^{\lambda Y}}{\E  e^{\lambda Y}}=0. 
\end{equation} 
By \cite{Petrov65}, the unique root of (\ref{equ:Legendre}) exists for any $x$ that satisfies $\E  Y<x<A_0=\lim_{\lambda\to B_-} \frac{\E  Ye^{\lambda Y}}{\E  e^{\lambda Y}}$ and $0<\lambda(x)<B$, where $B=\sup{\{\lambda: m(\lambda)<\infty\}}$.  Moreover, $\lambda(x)$ is a continuous and strictly monotone function of $x$, and therefore, $\lim_{x\to A_{0-}}\lambda(x)=B$, $\lim_{x\to {(\E  Y)}_+}\lambda(x)=0$, and $\lambda\bigl((\E  Y,A_0)\bigr)=(0,B)$. Differentiating $I^*(x)$ and using (\ref{equ:Legendre}), we have that
\begin{equation}
\label{equ:extrem}
\frac{d}{dx}I^*(x)=\lambda'(x)x+\lambda(x)-1-\frac{\lambda'(x)\E  Ye^{\lambda(x) Y}}{\E  e^{\lambda(x) Y}}=\lambda(x)-1.
\end{equation}
If $B\!>\!1$, there is a unique root $x_0$ of the equation $I^*(x)'\!=\!0$ such that $\lambda(x_0)\!=\!1$. Therefore, $\min_x I^*(x)\!=\!l(x_0,1)-x_0\!=\!0$. Otherwise, when $B=1$, 
\be
\rule{1ex}{0em}\inf\nolimits_x I^*(x)=\lim\nolimits_{x\to A_{0-}} \bigl(l(x,\lambda(x))-x\bigr)= \lim\nolimits_{x\to A_{0-}}\bigl(x\lambda(x)-\log m(\lambda(x))-x\bigr)=0.
\nnn
Note that $\frac{\E  Ye^{\lambda Y}}{\E  e^{\lambda Y}}=(\log m(\lam))'|_{\lam=1}$. Taking into account that 
$\log m(0)\!=\!\log m(1)\!=\!0$, and that $\log m(t)$ is a strictly convex function, we obtain 
 $\frac{\E  Ye^Y}{\E  e^Y}>0$. Then, the conclusion $x_0>0$ follows immediately from (\ref{equ:Legendre}), and the lemma is proved.  
\end{proof}

The following lemma allows to get a non precise asymptotic lower bound for $\E  e^{W_n}$ by using Doob's inequality.  

\begin{lemma}
\label{lem:lwbound}
Let the moment generating function $m(\lambda)=\E  e^{\lambda Y}$ exist in some neighbourhood of zero (the Cramer's condition), and suppose that $\Pr (Y>0)>0$. Then, for any $u>0$, there exists a sequence $\l k_n \r_1^{\infty}$ such that $k_n\topw{}{n\to\infty}\infty$, $k_n/n\topw{}{n\to\infty} 0$, and 
\begin{equation}
\label{equ:lwbound}
\liminf_{n\to\infty} e^{-k_n u}\E  e^{W_n}\geq 1-e^{-1}. 
\end{equation}
\end{lemma}

\begin{proof} 
First, we note that there $I(u)\in (0,\infty)$ for some $u>0$. Indeed, $\Pr (Y>0)>0$ implies that there exists $\delta>0$ such that 
$p_{\delta}=\Pr (Y>\delta)>0$. Then, from $\E  e^{\lambda Y}\geq p_{\delta} e^{\lambda\delta}$ for any $\lambda>0$, we conclude that 
\be
l(u,\lambda)=\lambda u-\log m(\lambda)\leq \lambda (u-\delta)-\log p_{\delta}.
\nnn 
Therefore, $\sup_{\lambda>0} l(u,\lambda)\leq -\log p_{\delta}<\infty$ for any $0<u\leq \delta$. Taking convexity of $\log m(\lambda)$ into account, with $\frac{\partial}{\partial \lambda}l(u,\lambda)|_{\lambda=0}=u-\E  Y>0$, we obtain that $\sup_{\lambda\leq 0}l(u,\lambda)=0$ for $u>0$. On the other hand,  conditions $\frac{\partial}{\partial \lambda}l(u,\lambda)|_{\lambda=0}>0$ and $l(u,0)=0$ imply that $I(u)>0$ for $0<u\leq \epsilon$ under some $\epsilon>0$. Hence, $I(u)\in (0,\infty)$ for $0<u\leq\delta\wedge\epsilon$. 

The Doob's inequality implies that
\begin{equation}\label{EexpWbyDoob}
\E  e^{W_n} \ge e^h\Pr \l \max_{[0,n]}(W_t) \ge h\r
\end{equation}
for any $h>0$.

To bound the probability $\Pr \l \max_{[0,n]}(W_t) \ge ku\r$, 
we again partition $[0,n]$ into $[n/k]$ subintervals of length $k$, for some integer $k$. According to the Cramer's large deviations theorem, the probability for the segment of random walk $S_k$ over any of these subintervals to exceed threshold $h=ku$ is
\[
\Pr (Y_1+\ldots+Y_k \ge ku ) =\exp\Bigl(-k\bigl(I(u)+\alpha(u,k)\bigr)\Bigr),
\]
where $\alpha(u,k)\topw{}{k\to\infty}0$ for any fixed $u>0$.  

Let $k_n=[x]$, where $x$ is a root of  the equation $xe^{x(I(u)+\alpha(u,x))}=n$ for some $u>0$ satisfying $I(u)\in (0,\infty)$. If at least one such segment of $S_n$ exceeds $h=ku$, it brings the maximum of $W_n$ at or above $h$. Since $k_n(I(u)+\alpha(k,u))\topw{}{n\to\infty}\infty$, we obtain that 
\be
\Pr (\max_{[0,n]}(W_t) \ge k_n u) &\!\ge\!& 1 - \bigl(1-e^{-k_n(I(u)+\alpha(u,k_n))}\bigr)^{[n/k_n]}
 \nn
&\!=\!&
1-\exp\left((-1+o(1))\frac{n}{k_n}e^{-k_n(I(u)+\alpha(u,k_n))}\right)\topw{}{n\to\infty} 1-e^{-1}.
\nnn
Then, for any $\epsilon>0$,
\be
 e^{-k_n u}  \E  e^{W_n}\geq\Pr (\max_{[0,n]}(W_t) \ge k_n u) \geq 1-e^{-1}-\epsilon,
\nnn
under a sufficiently large $n$, which implies (\ref{equ:lwbound}). 
\end{proof}

Result of Lemma~\ref{lem:lwbound} takes us to the following theorem.

\begin{theorem}
\label{prp:Cramers}
Assume the Cramer's condition for the distribution of $Y$ in some neighbourhood of zero, 
and suppose that $\Pr (Y>0)>0$. Then 
\begin{equation}
\label{equ:LBExp}
\liminf_{n\to\infty}  \frac{\E e^{W_n}}{n^{\gamma}}=\infty,
\end{equation}
for any $\gamma<\sup_{u>0} \l u/I(u)\r$.   
\end{theorem}

\begin{proof}  In  Lemma~\ref{lem:lwbound}, we proved that  (\ref{equ:lwbound})  holds for $k_n=[x]$, where $x$ is a root of the equation 
\be
xe^{x(I(u)+\alpha(u,x))}=n,
\nnn
for some $\alpha(u,x)\topw{}{x\to\infty}0$ and any fixed $u>0$ satisfying $I(u)\in (0,\infty)$. Since $x\topw{}{n\to\infty}\infty$, we conclude that 
\be
k_n=\log n /(I(u)+o(1))=\log n /I(u)+o(1), \ n\to\infty.
\nnn 
Hence,
\be
e^{-k_nh}=e^{-\bigl(\frac{u}{I(u)}+o(1)\bigr)\log n}=1/n^{\frac{u}{I(u)}+o(1)}.
\nnn
Taking into account the arbitrary choice of $u$, we obtain (\ref{equ:LBExp}) for any $\gamma<\sup_{u}\l u/I(u)\r$, which completes the proof.
\end{proof}

\begin{remark}(Uniqueness of the root). Extreme points of the function $I(u)/u$ solve the equation 
\begin{equation}
\label{equ:Ihextrem}
\left(\frac{I(u)}{u}\right)'=\left(\frac{I(u)-u}{u}\right)'=\frac{(I(u)-u)'-(I(u)/u-1)}{u}=0,
\end{equation}
if such solutions exist. Therefore, by (\ref{equ:extrem}), we obtain that 
$
\sup_u \l u/I(u)\r = 1/\lambda(u),
$
where $\lambda(u)=\lambda$ is a root of the equation 
\be
\frac{\E ( Ye^{\lambda Y} )}{\E  e^{\lambda Y}}=u. 
\nnn
Moreover,
\be
(I(u)/u)''=I''(u)/u+2(uI'(u)-I(u))/u^3,
\nnn
if $I$ is a twice differentiable function. Then there is at most one root $u_0$ of (\ref{equ:Ihextrem}) that is the point of minimum of $I(u)/u$, since $(I(u)/u)''|_{u=u_0}=I''(u_0)/u_0>0$. 
\end{remark}

By Lemma~\ref{lem:standard}, $\sup_{u>0} \l u/I(u) \r =1$ under $\E  e^Y=1$. Then, using Theorem~\ref{prp:Cramers}, we immediately obtain the following corollary.

\begin{corollary}
\label{cor:Cramers}
Under the conditions of Theorem~\ref{prp:Cramers}, (\ref{equ:LBExp}) holds for any $\gamma<1$.  
\end{corollary}

It is intuitively clear that the inequality in (\ref{equ:segmax}) is not optimal since $S_k\topw{\mbox{\scriptsize a.s.}}{k\to\infty}-\infty$. It should be more efficient using the segments of random length separated by successive zeroes of the CUSUM. Let $\ttau_0=0$ and  $\ttau_{k}=\inf\{t>\ttau_{k-1}:W_t=0\}$ for $k\in\Nb$, are the successive zeroes of the CUSUM process,  $M^{\circ}_k=\max_{[\ttau_{k-1},\ttau_k)} W_t$ and $l(n)=\max\{k:\ttau_k\leq n\}$. Then 
$$
\max_{[0,l(n)]} M^{\circ}_k\leq \max_{[0,n]} W_t \leq\max_{[0,l(n)+1]} M^{\circ}_k.
$$
Since $M^{\circ}_1,M^{\circ}_2,\ldots$ are i.i.d. random variables 
$$
\Pr(\max_{[0,m]} M^{\circ}_k>h)=1-(1-\Pr(M^{\circ}_1>h))^m
$$
for any fixed $m\in\Nb$. Using Fisher--Tippett--Gnedenko extreme value theorem and the well known convergence $l(n)/n\topw{\Pr}{n\to\infty}\mu^{-1}$ from the renewal theory, where $\mu=\E \tau_1$,  Iglehart \cite{Iglehart72} obtains the limit behaviour  of $\max_{t\leq n} W_t$ as $n\to\infty$ under the assumptions $\E e^{Y_1}=1$; $\E Y_1 e^{Y_1}<\infty$ and $Y_1$ is nonlattice. Theorem 2 in \cite{Iglehart72} implies that 
$$
\lim_{n\to\infty} \Pr_F(\max_{t\leq n} W_t-\log nb\leq x) =\Lambda^{1/\mu}(x), 
$$
where $\Lambda(x)=\exp(-e^{-x})$ is the standard Gumbel distribution function, $b$ is a normalizing constant and $\mu=\E \tau_1$, where $\tau_1$ is the first zero of the CUSUM process. Then 
\Be\label{IglehartApprox}
\lim_{n\to\infty} \Pr_F(\max_{t\leq n} W_t\leq x+\log nb) =\exp\Bigl(-\frac{1}{\mu}e^{-x}\Bigr), \; x\in\Rb. 
\Eeq
Hence, choosing $h=x+\log nb$ we obtain that for any $\lambda>0$
\begin{eqnarray*}
\frac{\E e^{\lambda W_n}}{n^{\lambda}}\geq b^{\lambda} e^{\lambda x}\Pr_F\bigl(\max_{t\leq n} W_t>x+\log nb\bigr)\rule{25ex}{0em}\\
\topw{\Pr}{n\to\infty} \exp\Bigl(\lambda(x+\log b)\Bigr)\Bigl(1-\exp\bigl(-\frac{1}{\mu}e^{-x}\bigr)\Bigr)
\end{eqnarray*}
Therefore, for any $\epsilon>0$,
\begin{equation}
\label{equ:lwbd}
\frac{\E e^{\lambda W_n}}{n^{\lambda}}\geq M-\epsilon\quad \mbox{under the sufficiently large}\quad n,
\end{equation}
where \vspace{-1.5em}
\begin{eqnarray*}
M&=&\sup_{x>0} \Bigl\{ \exp\Bigl(\lambda(x+\log b)\Bigr)\Bigl(1-\exp\bigl(-\frac{1}{\mu}e^{-x}\bigr)\Bigr)\Bigr\}\\
&=& \exp\Bigl(\lambda(x_0+\log b)\Bigr)\Bigl(1-\exp\bigl(-\frac{1}{\mu}e^{-x_0}\bigr)\Bigr);
\end{eqnarray*}
$x_0=-\log v_0\mu$ and $v_0$ is a root of the equation 
$
e^{-v_0}+v_0=\lambda.
$

In the particular case of $\lambda=1$ we obtain the linear lower bound for the exponential moment of CUSUM.

\section{Linear bounds for the exponential moments of CUSUM}
\label{sec:linbounds}

In this section, we apply Bell polynomials to study the asymptotic properties of the sequence of exponential moments of a CUSUM process. We start with the following lemma.  

\begin{lemma}
\label{lem:limSplus}
Let $Y_1,Y_2,\ldots$ be i.i.d. random variables satisfying \num{mgfcondition}. Then \hbox{$\E  e^{S^{+}_n}\!\topw{}{n\to\infty}\!2$.} 
\end{lemma}

\begin{proof}
The exponential moments of rectified sums can be expressed as,
\begin{eqnarray}
\label{equ:psm}
\E  e^{S^{+}_n}&=&\E  \Bigl(e^{S_n}\indi_{\{S_n\geq 0\}}+\indi_{\{S_n<0\}}\Bigr)\nonumber\\
&=&1+\E  (1-e^{S_n})\indi_{\{S_n<0\}}=1+\Pr (S_n<0)-\E  e^{S_n}\indi_{\{S_n<0\}}.
\end{eqnarray}
From \num{negEY}, we obtain convergence $S_n\topw{\Pr}{n\to\infty}-\infty$, and therefore, $\Pr(S_n<0)\topw{}{n\to\infty} 1$. This assures that the stopping time $\tau_{k}=\min\{n\in\Nb:\Pr(S_m<-k)\geq 1-1/k \;\mbox{for all}\;m\geq n\}$ is proper 
for every $k\in\Nb$. Let $a_n=k$ for all $n\in [\tau_k,\tau_{k+1})$. Then
\be
\E  e^{S_n}\indi_{\{S_n<0\}}&\leq & \E  \indi_{\{ -a_n\leq S_n<0\}}+e^{-a_n} \E  \indi_{\{S_n<-a_n\}}\!\nonumber \\
&\leq& \Pr (S_n\geq -a_n)+e^{-a_n}\leq 1/k+e^{-k}
\nnn
for all $n\geq \tau_k$. Convergence $\E  e^{S_n}\indi_{\{S_n<0\}}\topw{}{n\to\infty}0$ follows immediately. Combining the last argument with the convergence of $\Pr (S_n<0)\topw{}{n\to\infty} 1$, we obtain  that $\E  e^{S^{+}_n} \topw{}{n\to\infty} 2$ from~\num{equ:psm}.
\end{proof}

The following two lemmas state important properties of the sequence of rescaled Bell polynomials $\widetilde B_n$ defined in \num{rescaledBell} which we use to prove existence of a slant asymptote of~$\E e^{W_n}$. 

\begin{lemma}
\label{lem:Diff}
Let  $\Delta^{(k)}_n(x_1,\ldots,x_n)$ be the $k$-th order difference of $\widetilde B_n(x_1,\ldots,x_n)$, defined recursively as: 
\[
\Delta^{(k)}_0=B_0,\ \Delta^{(k)}_n(x_1,\ldots,x_n)=\Delta^{(k-1)}_n(x_1,\ldots,x_n)-\Delta^{(k-1)}_{n-1}(x_1,\ldots,x_n),
\]
for all $k,n\in\Nb$, with the convention $\Delta^{(0)}_n(x_1,\ldots,x_n)=\widetilde B_n(x_1,\ldots,x_n)$. Then, for all $r\in\Nb$,
\begin{equation}
\label{equ:BellDiff}
\Delta^{(r)}_n(x_1,\ldots,x_n)=\widetilde B(x_1-r,\ldots,x_n-r).
\end{equation}
\end{lemma}

\begin{proof}
To simplify notations, denote $b_n=\widetilde B(x_1,\ldots,x_n)$, $\Delta^{(1)}_{n}=b_n$, and $\Delta^{(k)}_{n}=\Delta^{(k)}_{n}(x_1,\ldots,x_n)$ for all $n\in\Nb\cup\{0\}$ and $k\in\Nb$.  
Then
\begin{eqnarray*}
\lefteqn{\Delta^{(1)}_{n+1} 
 = 
 b_{n+1}-b_n=\frac{1}{n+1}\sum_{k=0}^n b_{n-k}x_{k+1}-\frac{1}{n}\sum_{k=0}^{n-1} b_{n-k-1}x_{k+1}}
 \\
&=&\frac{b_0x_{n+1}}{n+1}+\frac{1}{n(n+1)}\Bigl\{n\sum_{k=0}^{n-1}(b_{n-k}-b_{n-k-1})x_{k+1}\Bigr\}-\frac{b_n}{n+1}\\
&=&\frac{b_0x_{n+1}}{n+1}+\frac{1}{n+1}\sum_{k=0}^{n-1}(b_{n-k}-b_{n-k-1})x_{k+1}-\frac{1}{n+1}\Bigl(b_0+\sum_{k=0}^{n-1}(b_{n-k}-b_{n-k-1})\Bigr)\\
&=&\frac{1}{n+1}\Bigl\{b_0(x_{n+1}-1)+\sum_{k=0}^{n-1}(b_{n-k}-b_{n-k-1})(x_{k+1}-1)\Bigr\}
 \nn
&=&
 \frac{1}{n+1}\sum_{k=0}^{n}\Delta^{(1)}_{n-k}(x_{k+1}-1) 
\end{eqnarray*}
Hence, 
\be
\sum_{k=0}^{n}\Bigl\{\Delta^{(1)}_{n+1}-\Delta^{(1)}_{n-k}(x_{k+1}-1)\Bigr\}=0.
\nnn
In a similar manner, we obtain by induction that for all $n,r\in\Nb$,
\be
\Delta^{(r)}_{n+1}=\frac{1}{n+1}\sum_{k=0}^{n}\Delta^{(r)}_{n-k}(x_{k+1}-r),
\nnn
and therefore, (\ref{equ:BellDiff}) holds.
\end{proof}

\begin{lemma}
\label{lem:Diff2}
Let $\{y_n\}_{n\in\Nb}$ be a numeric sequence with $\lim_{n\to\infty}y_n=0$. Then 
\be
\sum_{k=1}^{\infty} \widetilde B_k(y_1,\ldots,y_k)<\infty.
\nnn
\end{lemma}

\begin{proof}
Let $\widetilde y_n=\max_{k>n}|y_k|$ and $\undertilde{y}_n=\max_{k\leq n}|y_k|$ for all $n\in\Nb$. Then, for any $s\leq n\in\Nb$,
\begin{eqnarray*} 
\lefteqn{\bigl|\widetilde B_{n+1}(y_1,\ldots,y_n)\bigr|}  \nn
&&\leq  \frac{\undertilde y_s s}{n\!+\!1} \frac{1}{s}\sum_{k=0}^{s-1}|
 \widetilde B_{n-k}(y_1,\ldots,y_{n-k})|
 + \frac{\widetilde y_s(n\!-\!s\!+\!1)}{n+1} \frac{1}{n\!-\!s\!+\!1}\sum_{k=s}^{n}|\widetilde B_{n-k}(y_1,\ldots,y_{n-k})|\\
&&\leq \frac{\undertilde y_s s}{n+1}\max_{n-s+1\leq k \leq n} |\widetilde B_{k}(y_1,\ldots,y_k)|+\frac{\widetilde y_s(n-s+1)}{n+1}\max_{k\leq n-s} |\widetilde B_{k}(y_1,\ldots,y_{k})|\\
&&\leq \max_{k\leq n}|\widetilde B_k(y_1,\ldots,y_k)|\Bigl(\frac{s}{n+1}\,\undertilde y_s+\frac{(n-s+1)}{n+1}\widetilde y_s\Bigr)\\
&&=\max_{k\leq n}|\widetilde B_k(y_1,\ldots,y_k)|\Bigl(\widetilde y_s+\frac{s}{n\!+\!1}\,(\undertilde y_s-\widetilde y_s)\Bigr)\leq 
\max_{k\leq n}|\widetilde B_k(y_1,\ldots,y_k)|\Bigl(\widetilde y_s
 +\frac{s\undertilde y_{\infty}}{n\!+\!1}\Bigr).
\end{eqnarray*}
Now, choose $s_0\leq n_0$ so that $\widetilde y_{s_0}+\frac{s_0\undertilde y_{\infty}}{n_0+1}\leq 1/2$. Then, for any $n>n_0$,
\begin{equation}
\label{equ:d2bound}
 \widetilde B_n(y_1,\ldots,y_n)\leq \frac{1}{2^{n-n_0}}\max_{k\leq n_0} |\widetilde B_k(y_1,\ldots,y_k)|.
\end{equation}
Hence, 
\be
\sum_{k=1}^{\infty} \widetilde B_k(y_1,\ldots,y_k)\leq (n_0+2)\max_{k\leq n_0} |\widetilde B_k(y_1,\ldots,y_k)|<\infty,
\nnn
and the lemma is proved. 
\end{proof}
 
Using Lemmas \ref{lem:Diff} and \ref{lem:Diff2}, we obtain the following theorem.

\begin{theorem}
\label{pro:asymptote}
Let $\{x_n\}_{n\in\Nb}$ be a numeric sequence with the limit $\lim_{n\to\infty} x_n=2$ and let and $b_n=\widetilde B(x_1,\ldots,x_n)$ for all $n\in\Nb$. Then the sequence $\{b_n\}_{n\in\Nb}$ has 
an asymptote. 
\end{theorem}

\begin{proof} By (\ref{equ:BellDiff}), we obtain that $\Delta_n^{(2)}=\widetilde B_n(x_1\!-\!2,\ldots,x_n\!-\!2)$. By Lemma \ref{lem:Diff2}, there exists a finite limit
$
a=\lim_{n\to\infty} \Delta^{(1)}_n=\sum_{k=0}^{\infty} \Delta^{(2)}_k.
$
A slant asymptote exists if $a\not= 0$. 

Let $a=\lim_{n\to\infty} \Delta^{(1)}_n=\sum_{k=0}^{\infty} \Delta^{(2)}_k=0$ and $n_0\in\Nb$ is such that (\ref{equ:d2bound}) holds. Then, for any $s\in\Nb$, 
\be
\Delta^{(1)}_s+\sum_{k=s+1}^{\infty} \Delta^{(2)}_k=\sum_{k=0}^s \Delta^{(2)}_k+\sum_{k=s+1}^{\infty} \Delta^{(2)}_k=0,
\nnn
and therefore, $\Delta^{(1)}_s\!=\!-\sum_{k=s+1}^{\infty} \Delta^{(2)}_k$. The inequality $\Delta^{(2)}_k\!\leq\! \frac{\max_{s\leq n_0} |\Delta^{(2)}_s|}{2^{k-n_0}}$ for each $k\!\geq\! n_0$ implies that 
\be
\bigl|\Delta^{(1)}_k\bigr|\leq \sum_{k=s+1}^{\infty} \bigl|\Delta^{(2)}_k\bigr|\leq \frac{\max_{k\leq n_0} | \Delta^{(2)}_k|}{2^{k-n_0-1}}.
\nnn
Finally, 
\be
\Bigl|\sum_{k=0}^{\infty}\Delta^{(1)}_k\Bigr|\leq |b_{n_0}|+\sum_{k=n_0+1}^{\infty} \bigl|\Delta^{(1)}_k\bigr|\leq |b_{n_0}|+2\max_{k\leq n_0} | \Delta^{(2)}_k|<\infty. 
\nnn
Hence, the horizontal asymptote exists. 
\end{proof}

Note that $W_n\equiv 0$ if $\Pr(Y>0)=0$. In this case, $\E  e^{W_n}=1$, and the sequence of $\E  e^{W_n}$ has the horizontal asymptote $y=1$. In the next section we obtain that $\E  e^{W_n}$ has the slant asymptote in other cases, which follows immediately from Theorem~\ref{pro:asymptote}  and the unbounded property of this sequence. 

\section{General case and three types of behavior of CUSUM's exponential moments}\label{subsec:general}

Results of Sections~\ref{sec:Doob} and \ref{sec:linbounds} were obtained under the assumption \num{mgfcondition} 
on the exponential moment of original variables, $\E e^Y=1$. Here, we generalize these results to a broader class
of distributions. In this section, we replace condition \num{mgfcondition} on the common distribution of $Y_1,Y_2,\ldots$, by the following:
\Be\label{existenceofroot}
\mbox{Equation}\quad \E e^{\lambda Y}=1\quad \mbox{has a positive root}\quad \lamstar>0.
\Eeq
For example, for \num{existenceofroot} to hold, it is sufficient to have $\E Y < 0$ and $m(\lambda)=\E e^{\lambda Y}\to\infty$,
as $\lambda\to\infty$. Obviously, \num{mgfcondition} becomes a special case of \num{existenceofroot} with $\lamstar=1$.

To relate conditions \num{mgfcondition} and \num{existenceofroot}, introduce a sequence $Y_i^\ast=\lamstar Y_i$, which in turn generates a random walk $S_n^\ast = \lamstar S_n$ and a CUSUM process $ W_n^\ast = \lamstar W_n$. Condition \num{mgfcondition} is satisfied
by the sequence of $Y_i^\ast$, hence all the results of Sections~\ref{sec:Doob} and \ref{sec:linbounds} apply to the CUSUM $W_n^\ast$.

In particular, we have the lower and upper linear bounds for the exponential moments 
\[
 M_n^\ast(1)=\E e^{W_n^\ast} = \E e^{\lamstar W_n} = M_n(\lamstar),
\]
which establishes a connection with $M_n(\lam)$, the MGF of the original CUSUM process $W_n$. 

Actually, $\lamstar$ is the only value of $\lambda$ that leads to a linear behavior and a slant asymptote of $M(\lam)$, and it separates
three cases of the asymptotics of $M_n$.  
The following theorem characterizes three types of behavior of the CUSUM's moment generating function, as $n$ changes
between 0 and $\infty$.

\begin{theorem}
\label{thm:asymgf}
 Assume a sequence $\l Y_i\r_1^\infty$ of i.i.d. random variables with a common MGF $m(\lambda)=\E e^{\lambda Y}$ that satisfy conditions \num{existenceofroot}, $\E Y e^{\lambda Y}<\infty$, and $\omega=\lambda/\lambda^*>0$. Then, the MGF $M_n(\lambda)=\E e^{\lambda W_n}$  of  the induced CUSUM process $W_n$ has the following behavior in $n$:
 \be
 \mbox{For}\!&\!\lambda < \lamstar,\!&\! \mbox{ there exist $c_1,c_2>0$, such that $ c_1 n^\omega \leq M_n(\lambda)\le c_2 n^\omega$;} \nn
 \mbox{For}\!&\!\lambda = \lamstar,\!&\!M_n(\lambda) \mbox{ has a linear slant asymptote, as } n\to\infty; \nn
 \mbox{For}\!&\!\lambda > \lamstar,\!&\! \mbox{$\upsilon^n\le M_n(\lambda)\le c v^n$, where $c=v/(v-1)$ and $\upsilon = m(\lambda)>1$.} 
 \nnn 
 
 \begin{proof} 
The upper bound under $\lambda < \lamstar$ follows from the Lyapounov's inequality,
 \[
 \left( \E e^{\lambda W_n} \right)^{1/\lambda} =  \left( \E (e^{W_n})^\lambda \right)^{1/\lambda} 
 \le 
  \left( \E (e^{W_n})^{\lamstar} \right)^{1/\lamstar}  = \left( \E e^{\lamstar W_n} \right)^{1/\lamstar},
  \]
from where $M_n(\lambda) = \E e^{\lambda W_n} \le \left( \E e^{\lamstar W_n} \right)^{\lambda/\lamstar}
  = \left( M_n(\lamstar)\right)^{\lambda/\lamstar} \le (c^*n)^\omega=c_2 n^{\omega}$.

Theorem 2 in \cite{Iglehart72} states for non lattice distribution of $Y_1$ that 
$$
\lim_{n\to\infty} \Pr_F(\lambda^*\max_{t\leq n} W_t-\log nb\leq x) =\Lambda^{1/\mu}(x), 
$$
where $\Lambda(x)=\exp(-e^{-x})$ is the standard Gumbel distribution function, $b$ is a normalizing constant and $\mu=\E \tau_1$, where $\tau_1$ is the first zero of the CUSUM process. A~similar result for lattice distributions follows in the same manner from the asymptotic estimates in Feller (1971, vol. 2, sec. XI.6, page 377, (6.16)) with a version of the renewal theorem for lattice distributions (Feller, 1971, p. 360). Then 
$$
\lim_{n\to\infty} \Pr_F\Bigl(\lambda^*\max_{t\leq n} W_t\leq x+\log nb\Bigr) =\exp\Bigl(-\frac{1}{\mu}e^{-x}\Bigr), \; x\in\Rb. 
$$
Hence, applying the first inequality in (\ref{equ:GenDoob}) under the fixed $h=\frac{x}{\lambda^*}+\frac{1}{\lambda^*}\log nb$ we obtain that for any $\lambda>0$
\begin{eqnarray*}
\frac{\E e^{\lambda W_n}}{n^{\omega}}\geq b^{\omega} e^{\omega x}\Pr_F\bigl(\lambda^*\max_{t\leq n} W_t>x+\log nb\bigr)\rule{25ex}{0em}\\
\topw{\Pr}{n\to\infty} \exp\Bigl(\omega(x+\log b)\Bigr)\Bigl(1-\exp\bigl(-\frac{1}{\mu}e^{-x}\bigr)\Bigr).
\end{eqnarray*}
Therefore, for any $\epsilon>0$,
\begin{equation}
\label{equ:lwbd}
\frac{\E e^{\lambda W_n}}{n^{\omega}}\geq M-\epsilon\quad \mbox{under the sufficiently large}\quad n,
\end{equation}
where \vspace{-1.5em}
\begin{eqnarray}
M&=&\sup_{x>0} \Bigl\{ \exp\Bigl(\omega(x+\log b)\Bigr)\Bigl(1-\exp\bigl(-\frac{1}{\mu}e^{-x}\bigr)\Bigr)\Bigr\}\label{constM}\\
&=& \exp\Bigl(\omega(x_0+\log b)\Bigr)\Bigl(1-\exp\bigl(-\frac{1}{\mu}e^{-x_0}\bigr)\Bigr);
\nonumber\end{eqnarray}
$x_0=-\log v_0\mu$ and $v_0$ is a root of the equation 
$
e^{-v_0}+v_0=\omega.
$
The upper bound under $\lambda < \lamstar$ holds. 

Existence of a slant asymptote for $M_n(\lambda^*)$ as $n\to\infty$ follows immediately from Theorem~\ref{pro:asymptote}  and (\ref{equ:lwbd}) under $\lambda=\lambda^*$.

From 
$
\widetilde{B}_n(x,\ldots,x^n)=x^n
$
we obtain by Lemma \ref{lem:Diff} that for any $x>0$ 
$$
\widetilde{B}_n(x+1,\ldots,x^n+1)=\sum\nolimits_{k=0}^n x^k.
$$
Under $\lambda>\lambda^*$  similarly to (\ref{equ:psm}) we conclude that
$$
\E e^{\lambda S_n^+} \leq 1+\upsilon^n.
$$
Then taking into account the increasing property of Bell polynomials we obtain the upper  bound
$$
M_n(\lambda)=\widetilde B_n(\E e^{\lambda S^+_1},\ldots,\E e^{\lambda S^+_n})\leq \widetilde B_n\bigl(1+\upsilon,\ldots,1+\upsilon^n\bigr)\leq \sum\nolimits_{k=0}^n \upsilon^k\leq v^n \frac{v}{v-1}.
$$
Finally, lower bound under the case of $\lambda > \lamstar$ follows from a simple observation that $W_n\ge S_n$ for all $n$.  Therefore,
  \[
  M_n(\lambda) \ge \E e^{\lambda S_n} = (\E e^{\lambda Y})^n = \upsilon^n,
  \]
where $\upsilon=m(\lambda)>1$ for $\lambda > \lamstar$ by the convexity of $m(\lambda)$.
 \end{proof}
 
\end{theorem}

\begin{figure}
\begin{picture}(0,130)(240,0)
\put(33,-18){\includegraphics[scale=.39]{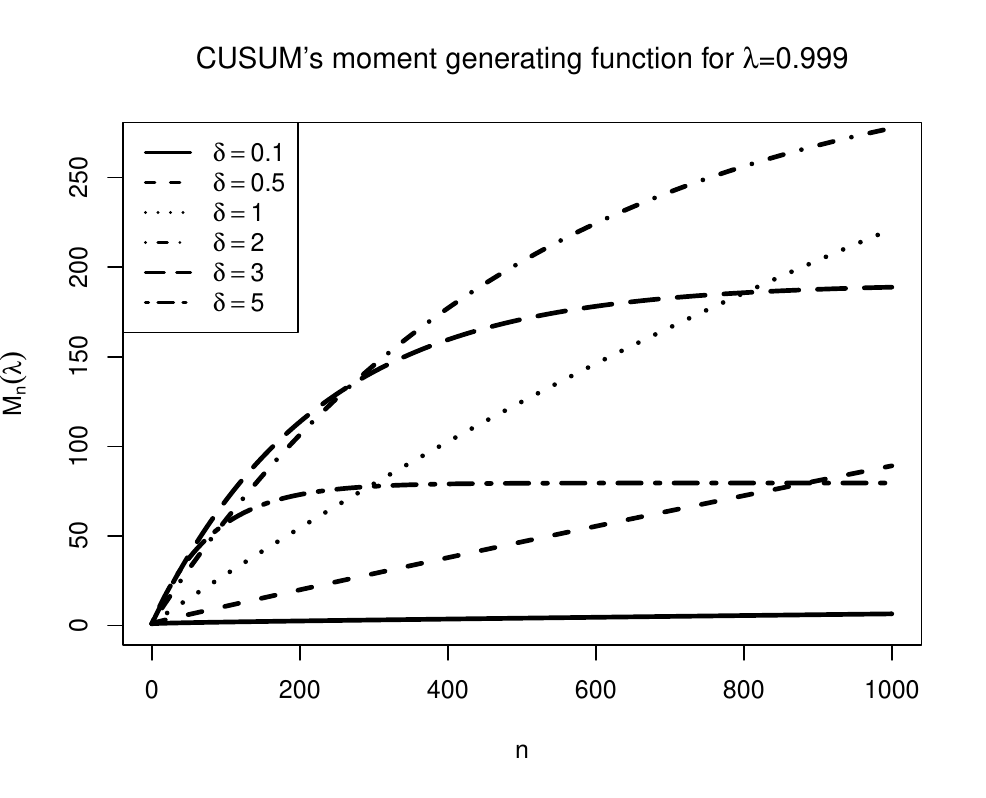} }
\put(247,-18){\includegraphics[scale=.39]{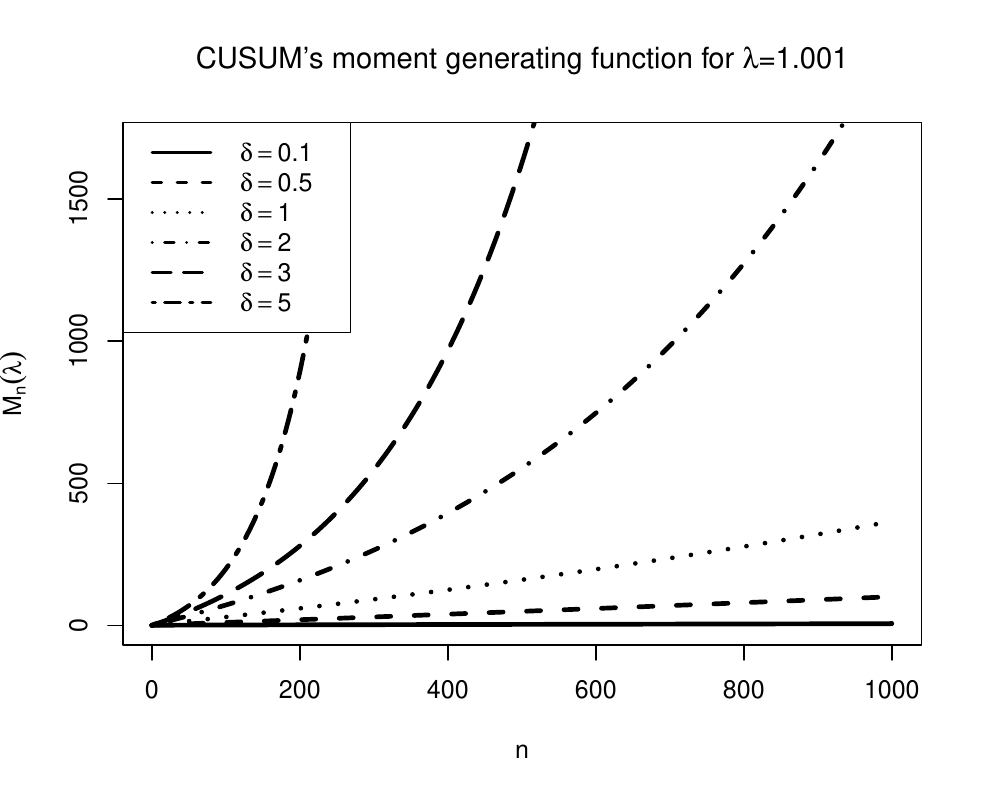}}
\end{picture}
\caption{CUSUM's moment generating functions for $\lambda < \lamstar$ (left) and for
$\lambda > \lamstar$ (right).\label{fig:mgfothertwo}}
\end{figure}

We have already observed the linearity of CUSUM's exponential moments corresponding to the case $\lambda=\lamstar=1$ in
Figure~\ref{fig:expmoments}. The other two types of behavior are seen in Figure~\ref{fig:mgfothertwo}. Note the evident
slow growth for $\lambda = 0.999 < \lamstar$ and exponential growth for $\lambda = 1.001 > \lamstar$ even for
a very slight deviation from the borderline quantity $\lamstar$. As in Figure~\ref{fig:expmoments}, $\delta$ on these graphs
represents the detectable difference between two distributions in a change-point problem.

\section{Practical applications}\label{sec:appl}
Obtained results and computational methods have applications in several practical
fields.

\subsection{Testing for abrupt changes}\label{subsec:chp}
A classical change-point problem refers to situations where a distribution $F$ of
observed data is abruptly replaced by a distribution $G$ at an unknown moment $\nu$, 
\[
\l\ba{lll}
 X_1,\ldots,X_\nu & \sim & F \\
 X_{\nu+1},\ldots,X_n & \sim & G 
 \ea\right.
\]
The likelihood-ratio test for occurrence of a change-point,
\[
H_0: \, X_1,\ldots,X_n \sim F \mbox{ against } 
H_A: \,  X_1,\ldots,X_\nu \sim F; \
 X_{\nu+1},\ldots,X_n \sim G \mbox{ for some } \nu ,
\] 
is based on the CUSUM statistic 
\be\label{chp_statistic}   
W_n
&=&
 \log\frac{\d\max \mbox{\eufm{L}}(X_1,\ldots,X_n\c\ H_0\cup H_A)}{\mbox{\eufm{L}}(X_1,\ldots,X_n\c\ H_0)}  
 \nn
&=&
\log\frac{\d\max_{\nu\in[0,n]}\left(\prod_1^\nu f(X_i)\prod_{\nu+1}^n g(X_i)\right)}{\d\prod_1^n f(X_i)} =
 \max_{\nu\in[0,n]} \sum_{\nu+1}^n \log \frac{g(X_i)}{f(X_i)}.
\eeq
It matches the definition \num{equ:CUSUM} of the CUSUM process $W_n$, with increments
$Y_i = \log \frac{g(X_i)}{f(X_i)}$ representing log-likelihood ratios of densities $f$ and $g$
of the default and the disturbed distributions $F$ and $G$ with respect to some 
reference measure $\mu$. Introduced in \cite{Page54}, this test has become standard for
testing occurrence of abrupt changes (\cite{BassevilleNikiforov93}, sec. 2.2; 
\cite{TartakovskyNikiforovBasseville14}, sec. 8.2; and many others).

The test rejects $H_0$ in favor of $H_A$ and declares detection of a change-point 
for large values of $W_n \ge h$. Then, controlling the probability of Type I error,
which would result in a false alarm, at level $\alpha$ reduces to a problem of finding a threshold $h$
such that
\Be\label{TypeIcontrol}
\P\l W_n \ge h\c\ H_0\r = \P_F\l W_n \ge h\r \le \alpha.
\Eeq
This can be achieved by bounding the probability in \num{TypeIcontrol} by means
of Markov, Chebyshev, or Chernoff inequality and using, respectively, the moments 
\num{meanW}, \num{varianceW}, or \num{expW} that are explicitly established in Theorem~\ref{thm:moments}
or calculated by the recursive or matrix method proposed in Sec.~\ref{sec:recursive}.

For example, we can bound this probability as
\Be\label{boundTypeI}
\P\l W_n\ge h \r \le e^{-h} \E e^{W_n} = e^{-h}M_n(1),
\Eeq
Bounding the right-hand side of \num{boundTypeI} by $\alpha$ and solving for $h$, we obtain the 
threshold 
\Be\label{threshold_EM}
h(\alpha) = \log\frac{M_n(1)}{\alpha}
\Eeq
that is guaranteed to control the false alarm probability within the level $\alpha$.
The exponential moment $M_n(1)$ can be calculated by methods \num{exponentialWrecursiveThm} or \num{comp_matrix}. A recursive computation of $M_n(\lam)$ for Normal distributions is detailed in Sec.~\ref{subsec:expmoments}.

\subsection{Testing for transient changes}

Data sequence $\l X_k\r$ experiences a {\sl transient change} if its initial distribution
$F$ changes to another distribution $G$ and then returns to the original distribution,
\[
\l\ba{lll}
 X_1,\ldots,X_\alph & \sim & F \\
 X_{\alph+1},\ldots,X_\bet & \sim & G \\
 X_{\bet+1},\ldots,X_n & \sim & F
 \ea\right.
\]
A likelihood-ratio test for the occurrence of a transient change,
$H_0:  X_1,\ldots,X_n \sim F$ against   
$H_A:   X_1,\ldots,X_\alph$, $X_{\bet+1},\ldots,X_n \sim F$; 
$X_{\alph+1},\ldots,X_\bet \sim G$ for some $\alph$ and $\bet$, 
is based on the 
{\sl maximum-CUSUM} test statistic 
\cite{BaronMalov23transient, levin1985cusum, siegmund1986boundary}
\be   
\Lambda_n 
&=& \log\frac{\d\max \mbox{\eufm{L}}(X_1,\ldots,X_n\c H_0\cup H_A)}{\mbox{\eufm{L}}(X_1,\ldots,X_n\c H_0)} 
= \log\frac{\d\max_{0\le\alph<\bet\le n}\!\left(\prod_1^\alph f(X_i)\prod_{\alph+1}^\bet g(X_i)\prod_{\bet+1}^n f(X_i)\right)}{\d\prod_1^n f(X_i)} 
 \nn[-5mm]
&=&\max_{0\le\alph<\bet\le n} \sum_{\alph+1}^\bet \log\frac{g(X_i)}{f(X_i)}
= \max_{0\le\alph<\bet\le n} (S_\bet - S_\alph)
= \max_{0<\bet\le n} (S_\bet - \min_{0\le\alph<\bet}S_\alph)
= \max_{0<\bet\le n} W_\bet.
\nnn

We have studied the maximum-CUSUM statistic $\max_{[0,n]}W_k$ 
in Sections~\ref{sec:Doob} and \ref{sec:linbounds}.
All the assumptions made in these sections are satisfied under the default distribution
$F$ in the abrupt and transient change-point situations:
\[
\E_F(Y_i) = \E_F\log\frac{g(X_i)}{f(X_i)} < 0,
\]
by the Kullback inequality \cite[(4.2)]{Kullback54}, and
\[
\E_F\left( e^{Y_i}\right) = \E_F \frac{g(X_i)}{f(X_i)} = \int g(x)d\mu(x) = 1.
\]
Then, all the results of Sections~\ref{sec:Doob} and \ref{sec:linbounds} apply,
$U_t = e^{ W_t }$ is a {submartingale}, and the maximal inequality 
\num{equ:DoobCUSUM} holds, relating the probability of a false alarm with
the CUSUM's exponential moment:
\[
\P\l \mbox{false alarm}\r = \P_F\l \max_{[0,n]}W_k \ge h\r \le e^{-h}\E_F(U_n).
\]
Hence, threshold \num{threshold_EM} controls the false alarm probability for the transient
changes as well!

We can offer a further bound for the probability of a false alarm, leading to a universal threshold.
Here we make use of identities for the incomplete exponential Bell polynomials,
\[
B_{nk}(cx_1,\ldots,cx_{n-k+1})=c^k B_{nk}(x_1,\ldots,x_{n-k+1}) \mbox{ for any } c>0,
\] 
and 
\[
B_{nk}(0!,\ldots,(n-k)!)=|s(n,k)|=\left[n\atop k\right]
\] 
(see \cite{Comtet74}, Theorem B(3i)), 
where $s(n,k)$ are the Stirling numbers of the first kind. 
From (\ref{equ:bellexp}), we obtain that 
\be
\E_F e^{W_n}\leq \frac{1}{n!} \sum_{k=1}^n \mbox{$\left[n\atop k\right]$} 2^k=\frac{2\cdot 3\cdot\ldots\cdot (2+n-1)}{n!}=n+1,
\nnn
and, therefore, 
\Be\label{bound1}
\P\l \mbox{false alarm}\r \le e^{-h}\E_F(U_n) \leq e^{-h}(n+1). 
\Eeq
Thus, selecting a threshold
 \Be\label{universalUB}
 h(\alpha) = \log\frac{n+1}{\alpha}
 \Eeq
guarantees control of the false alarm probability in transient change-point detection
within the level $\alpha$,
$\P\l \mbox{false alarm}\r \le \alpha$.
This is a universal threshold, it serves all the distributions and tops more accurate thresholds 
of the form $h(\alpha) = \log(\alpha^{-1}M_n(1))$.

\subsection{Tighter change-point detection threshold}

Application of the Doob–Meyer decomposition to the change-point scenario 
allows to reduce the threshold $h$ further, 
without losing control of the false alarm probability. 

Evaluating the right-hand side of \num{equ:marbound} for our change-point scenario with
$Y = \log \frac{g}{f}(X)$, we obtain
\be 
\E_F \left(1-e^Y\right)^+ 
&=&
 \E_F \left(1-\frac{g(X)}{f(X)}\right)\indi_{\l f(X) > g(X)\r}
 \nn
&=&
 \P_F\l f(X) > g(X)\r -  \P_G\l f(X) > g(X)\r = D_{F,G}.
\nnn
It yields an upper bound for the probability of a false alarm,
\Be\label{PFA_prob_difference}
\Pr\l \mbox{false alarm}\r 
= \Pr_F\{ \max_{[0,n]} W_t \ge h \} 
\leq e^{-h}\bigl(1+nD_{F,G}\bigr), 
\Eeq
which is a tighter bound than \num{bound1}.

We conclude that for any $\alpha\in (0,1)$, the threshold 
 \Be\label{h_D}
 h_\alpha = \log\frac{1+nD_{F,G}}{\alpha}
 \Eeq
guarantees the control of the probability of a false alarm at level $\alpha$,
\be
\Pr\l \mbox{false alarm}\r = \Pr_F \bigl(\max_{0\le t\le n} W_t \ge h_{\alpha}\bigr) \leq \alpha.
\nnn

The quantity $D_{F,G} = \P_F\l f(X) > g(X)\r -  \P_G\l f(X) > g(X)\r$, appearing in the 
probability of a false alarm \num{PFA_prob_difference}, represents a measure of discrepancy between
the distributions $F$ and $G$. It is non-negative, because
\begin{eqnarray}
\P_F\l f(X) \!>\! g(X)\r&=& \int_{f(x) > g(x)}\! f(x) d\mu(x)\nonumber \\
&\ge & \int_{f(x) > g(x)}\! g(x) d\mu(x) = \P_G\l f(X) \!>\! g(X)\r,
\label{positiveD}
\end{eqnarray}
and $D_{F,G}=0$ when \num{positiveD} is an equality, which occurs only when $\P_F\l f(X)>g(X)\r=0$, or 
$f\equiv g$, $\mu$-almost everywhere. Small changes, with a low $D_{F,G}$, are generally harder to detect,
the CUSUM process $W_n$ is less steep, and we see that such changes should be detected with lower thresholds
given by \num{h_D}.

\subsection{Change-point detection in randomly stopped sequences}

The upper bounds in (\ref{equ:marbound}) and \num{PFA_prob_difference} can be extended from a fixed $n$ to 
finite Markov stopping times $T$ with respect to filtration $\Ff$. 
Applying the Lenglart's inequality \cite[Theorem 1]{Lenglart77},
\Be\label{PFA_Markov}
\Pr\Bigl(\max_{[0,T]} W_t \ge h\Bigr)\leq  \frac{\E \widetilde U_T}{e^h}
 \leq  \frac{1+D_{F,G}\E_F T }{e^h}.
\Eeq

A practical example is a situation when the manufacturing process is stopped because of an observed extreme measurement,
such as $T = \inf\l n: |X_n| \ge \mu_X + 3\sigma_X\r$, 
or when an alarm is triggered by 
any of the Western Electric, AIAG, Nelson, Boeing, or General Electric control chart rules that are popular in 
statistical quality control \cite{western1958statistical, Montgomery91,
nelson1984shewhart, wheeler1995advanced}.
In this situation, the main question is whether the extreme was caused by a structural change in the process or by a
random fluctuation that is consistent with the default distribution. The test statistic $W_T$ given by \num{chp_statistic} 
can then be used with $n=T$, and its probability of a Type I error is bounded by \num{PFA_Markov}.

\subsection{Sequential detection and comparison of stopping times}

Detecting a change-point sequentially, one observes a data sequence $X_1,X_2,\ldots$ and stops
at the first moment when the CUSUM process exceeds threshold $h$,
\[
T_h = \inf\l n:\,W_n\ge h\r.
\]
This stopping rule is well studied, and it has a number of optimal properties, see 
\cite{Page54,Lorden71,Moustakides86,Ritov90}.
Result \num{PFA_Markov} allows comparison of $T_h$ with other stopping rules. For example,
the probability of $T_h$ occurring during any given ``excursion'' of a CUSUM, the segment between successive visits to 0, is 
bounded by
\[
\Pr\l T_h < \ttau_1 \r = \Pr\Bigl(\max_{0\le t\le \ttau_1} W_t \ge h\Bigr)\leq  \frac{1+D_{F,G}\E_F \ttau_1 }{e^h},
\]
and the probability of $T_h$ prior to the $k$-th visit to 0 is
\[
\Pr\l T_h < \ttau_k \r = \Pr\Bigl(\max_{0\le t\le \ttau_k} W_t \ge h\Bigr)\leq  \frac{1+D_{F,G}\E_F \ttau_k }{e^h}
= \frac{1+D_{F,G}k\E_F \ttau_1 }{e^h},
\]
where $\ttau_0=0$ and  $\ttau_{k+1}=\inf\{t>\ttau_k:W_t=0\}$ are the successive zeroes of the CUSUM process.

\subsection{Comparison of thresholds}\label{subsec:comparison}

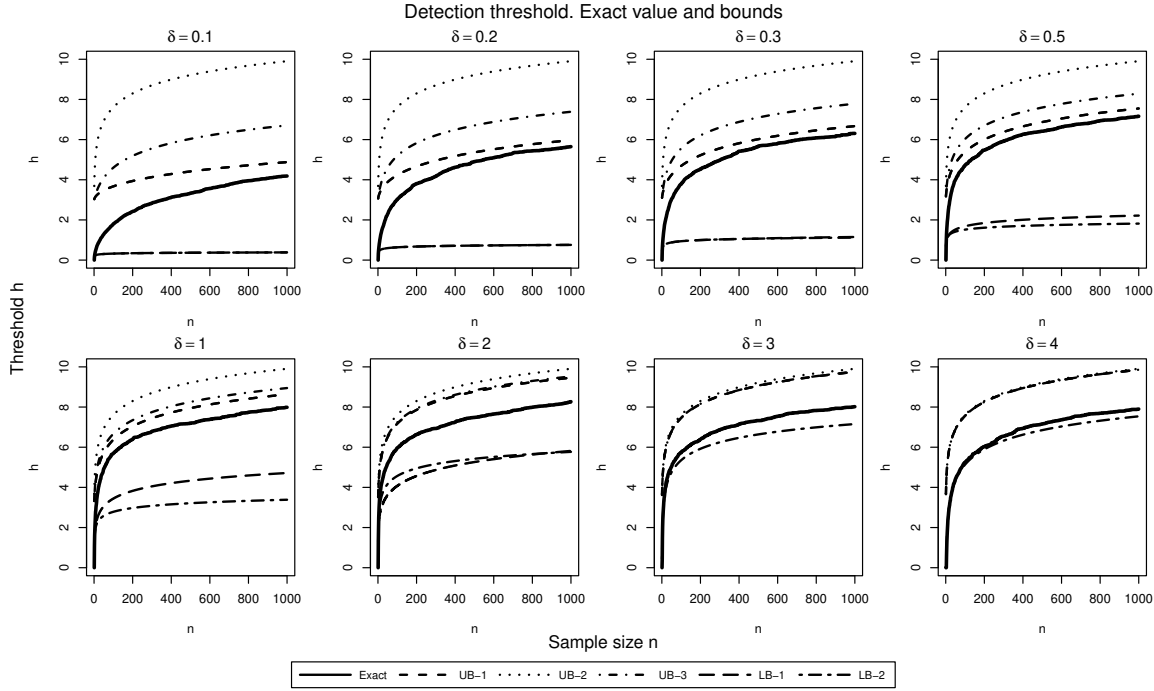
\begin{figure}\scalebox{0.92}[0.90]{
\begin{picture}(0,150)(230,-20)
\put(-10,0){\scalebox{.25}{\includegraphics{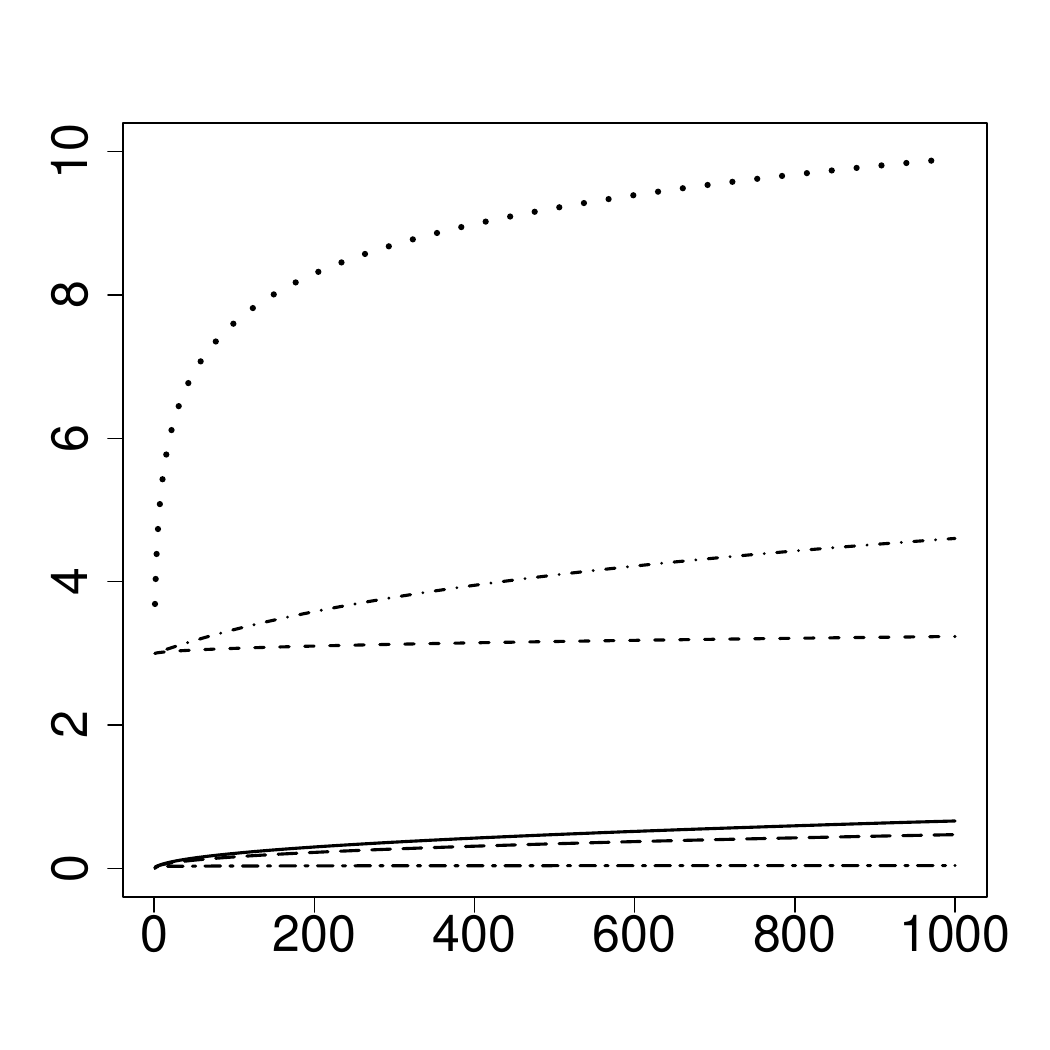}}}
\put(110,0){\scalebox{.25}{\includegraphics{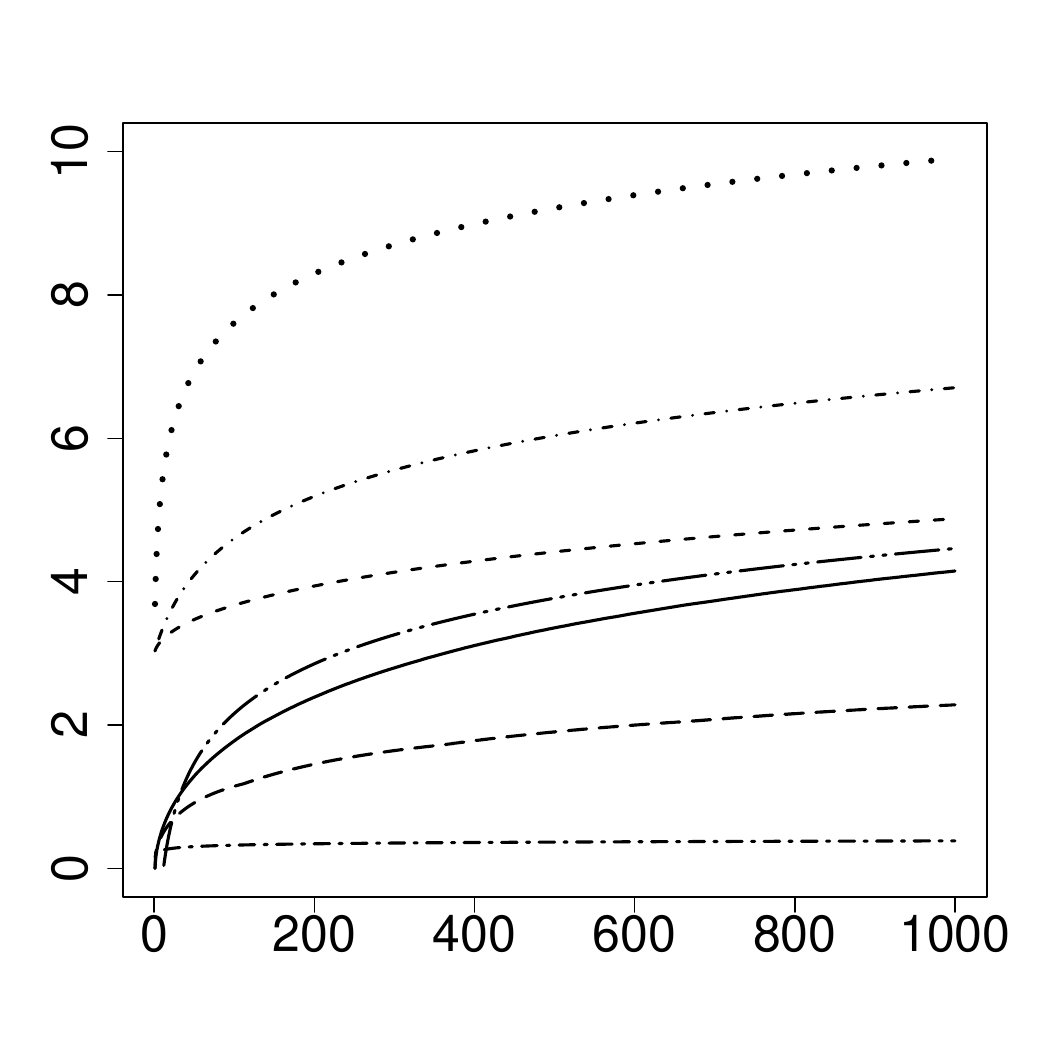}}}
\put(230,0){\scalebox{.25}{\includegraphics{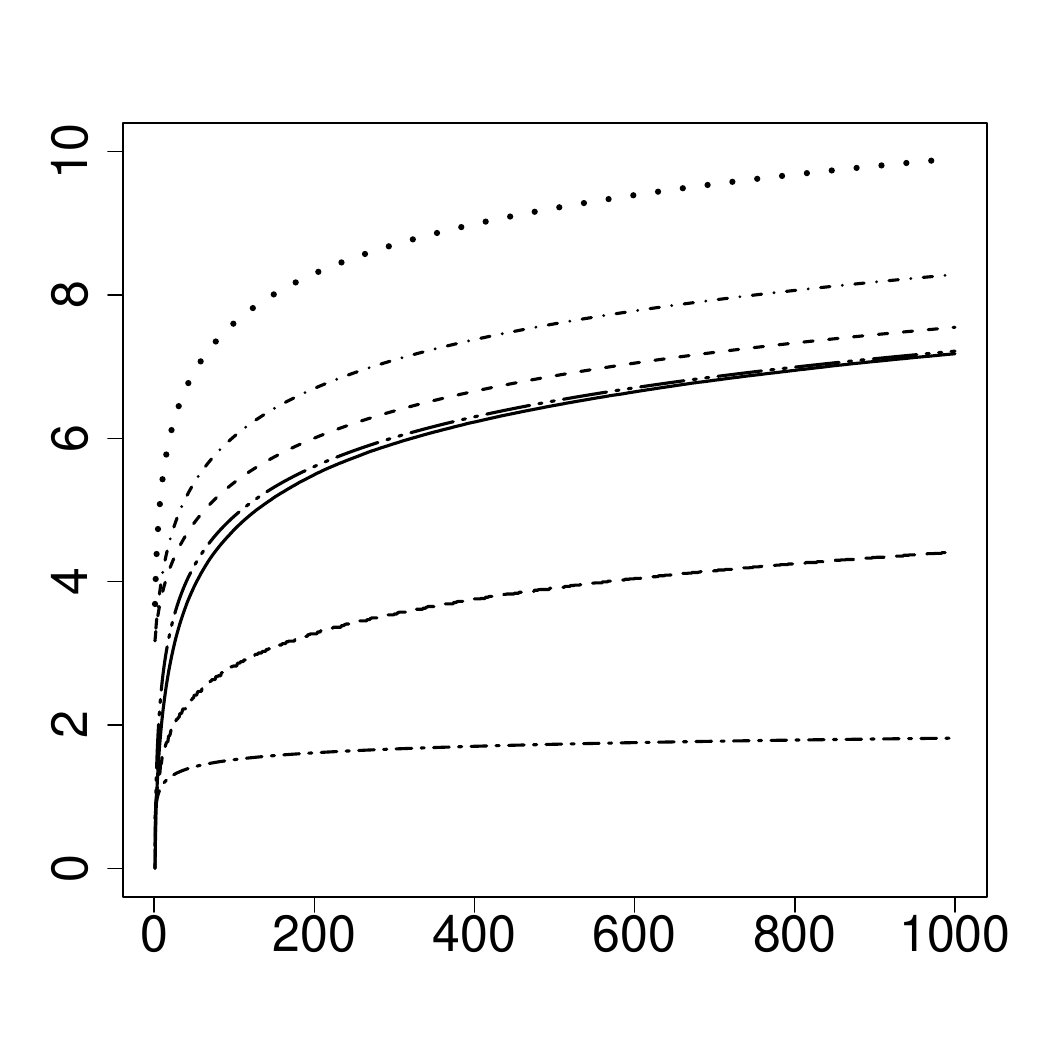}}}
\put(350,0){\scalebox{.25}{\includegraphics{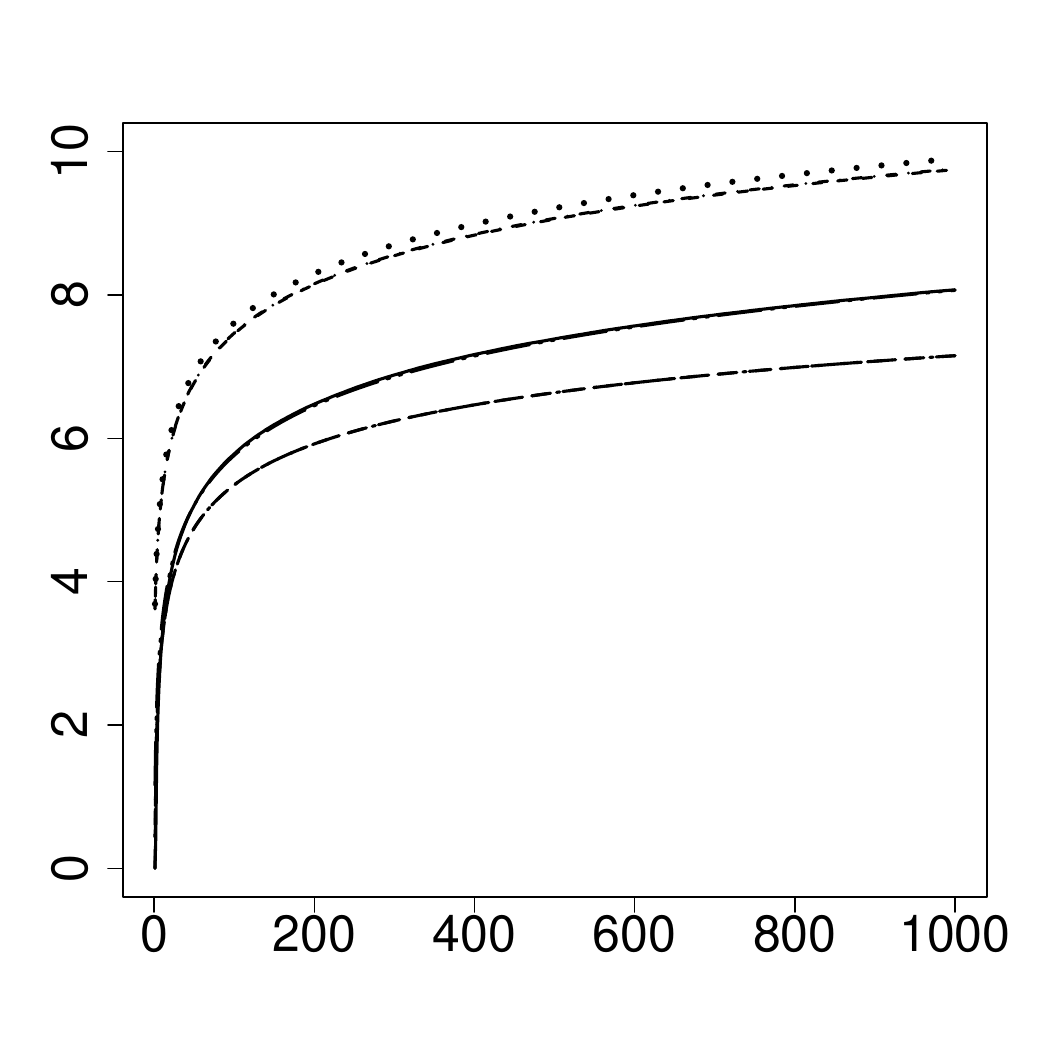}}}
\put(-10,35){\begin{rotate}{90}{\footnotesize Threshold} $h$\end{rotate}}
\multiput(25,0)(120,0){4}{{\footnotesize Sample size} $n$}
\put(37,117){$\delta = 0.01$}
\put(159,117){$\delta = 0.1$}
\put(279,117){$\delta = 0.5$}
\put(399,117){$\delta = 3.0$}
\put(-5,-24){\begin{tikzpicture}

\draw[solid, thick] (0,0) -- (1,0);
\node[right] at (1.1,0) {\small Exact};

\draw[dashed, thick] (2.5,0) -- (3.5,0);
\node[right] at (3.6,0) {\small UB-1};

\draw[dotted, thick] (5,0) -- (6,0);
\node[right] at (6.1,0) {\small UB-2};

\draw[dash dot, thick] (7.5,0) -- (8.5,0);
\node[right] at (8.6,0) {\small UB-3};

\draw[line width=0.75pt,dash pattern=on 6pt off 2pt on 2pt off 2pt]
      (10,0) -- (11,0);
\node[right] at (11.1,0) {\small LB-1};

\draw[line width=0.75pt,dash pattern=on 9pt off 3pt]
      (12.5,0) -- (13.5,0);
\node[right] at (13.6,0) {\small LB-k};

\draw[line width=0.75pt,dash pattern=on 10pt off 2pt on 1.5pt off 2pt on 1.5pt off 2pt]
      (15,0) -- (16,0);
\node[right] at (16.1,0) {\small IT};

\end{tikzpicture}}

\end{picture}}
\caption{Thresholds for the transient change-point detection and their proposed upper and lower bounds.
Four scenarios correspond to four detectable differences $\delta$ from 0.01 to 3.0. Thresholds are calculated
to control the false alarm rate within $\alpha=0.05$. The bounds
are explained in Section~\ref{subsec:comparison}.}
\label{fig:threshold_bounds}
\end{figure}

Detection thresholds for the transient change-point and the maximum-CUSUM statistic 
are depicted in Figure~\ref{fig:threshold_bounds}, along with 
our developed upper and lower bounds.


For each detectable difference $\delta$, the solid line on Figure~\ref{fig:threshold_bounds}
represents the exact threshold computed by simulating upper $\alpha$-quantiles of the distribution
of the maximum-CUSUM statistic $\max_{[0,n]} W_t$. It is shown along with the Doob's maximal inequality induced
upper bound UB-1, given in
\num{threshold_EM}; the universal upper bound UB-2, obtained from Bell's polynomials and given in \num{universalUB}; 
and its improvement UB-3, given in \num{h_D} and showing a discrepancy measure $D_{F,G}$. 
The lower bound
LB-k is calculated from \num{equ:segmax}. For Normal distributions, inequality \num{equ:segmax} results in
\Be\label{ineq_norm}
\Pr\l\max_{[0,n]}W_t\ge h\r \ge 1 - \Phi^{[n/k]}\left(\frac{h+k\delta^2/2}{\delta\sqrt{k}}\right),
\Eeq
leading to the threshold 
\Be\label{h-LB1}
h \ge \delta\sqrt{k}\,\Phi^{-1}\left((1-\alpha)^{1/[n/k]}\right) - k\frac{\delta^2}{2}
\Eeq
that controls the false alarm probability at level $\alpha$. Maximization of \num{h-LB1} over integer $k$ yields the lower bound LB-k.
Its simplified version, with $k=1$, results in the lower bound LB-1, following (22) in \cite[Sec. 5.2]{BaronMalov23transient}. 
Finally, the approximate upper bound IT is based on the Iglehart's asymptotic result \num{IglehartApprox}.

Figure~\ref{fig:threshold_bounds} offers several insights. First, detection thresholds are dependent on both the
detectable difference $\delta $ and the sample size $n$; they are increasing functions of $n$. Second, upper bounds appear more accurate
for {moderate} differences $\delta$ while lower bounds are { quite sharp} for {small} $\delta$. To maintain the control of 
a false alarm rate, only the threshold's upper bounds can be used. Also, small and moderate differences represent the
most meaningful change-point problems, because differences of 2 or more standard deviations are detectable with a naked eye.
Third, we can see that the upper bounds are compatible with the lowest thresholds that maintain the false alarm probability
at a given level. Within the mentioned probability constraint, low thresholds are preferred in order to minimize the detection
delay.
Although thresholds such as \num{threshold_EM} are based on the Doob's maximal inequality, they appear nearly optimal, for moderate $\delta$ and larger $n$, {whereas under larger $\delta$ they are closing to the universal bound}. {
The approximate upper bound IT is not displayed in Figure~\ref{fig:threshold_bounds} for $\delta=0.01$ because it takes negative values when $n<1038$, making it inapplicable in this case. In the remaining cases, the asymptotic threshold IT shown in Figure~\ref{fig:threshold_bounds} provides a close approximation to the exact threshold, but under an extremely large $\delta$ this approximation significantly overestimates the theoretical threshold. 
}

\begin{figure}
\begin{picture}(0,150)(240,0)
\put(80,150){\begin{rotate}{270}\scalebox{.3}[.39]{\includegraphics[scale=1]{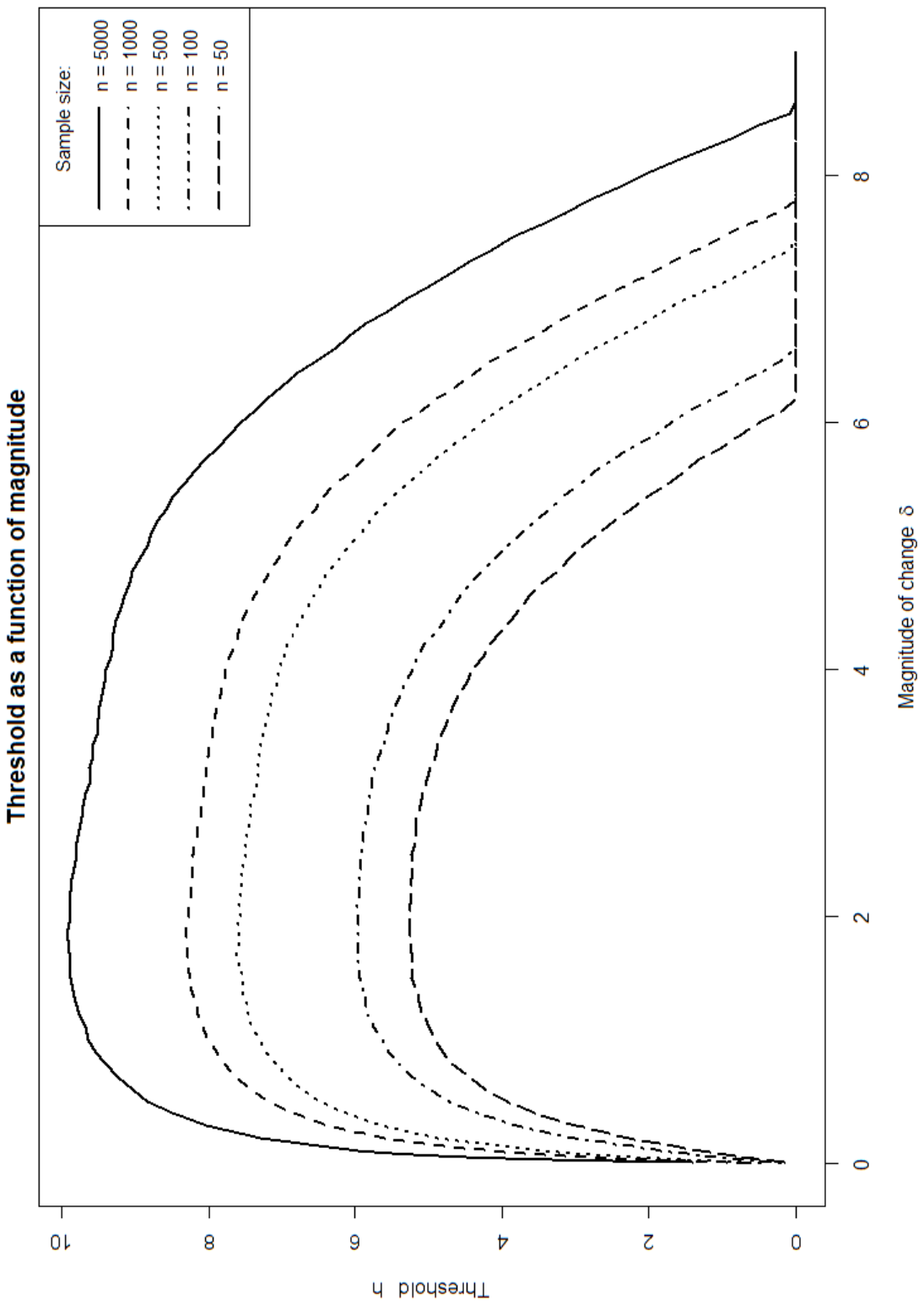} }
\end{rotate}}
\end{picture}
\caption{Threshold $h$ as a function of the detectable difference $\delta$, for various sample sizes $n$ \label{fig:thresholds_of_delta}}
\end{figure}

It is also informative to graph detection thresholds as functions of $\delta$, the magnitude of change, Figure~\ref{fig:thresholds_of_delta}.
Thresholds are low for very small as well as very large $\delta$.
Indeed, small $\delta$ cause a small negative mean $-\delta^2/2$ and a small standard deviation $\delta$ of increments $Y_t$. This results in a rather flat random walk $\l S_t\r$ and the induced CUSUM process $\l W_t\r$, with small values and small quantiles. On the other hand, large $\delta$, representing easy detection cases, imply large negative increments $Y_t$, sending the CUSUM process $W_t$ to 0 for most $t$, which again leads to low quantiles of $\max_{[0,n]}W_t$.

\subsection{Detection thresholds for unbounded sequences}\label{sec:mean_delay}

Thresholds (\ref{threshold_EM}, \ref{universalUB}, \ref{h_D}, \ref{h-LB1}) 
control the probability of a false alarm, 
$\P\l T_h\le \nu\r$, during the time segment $[0,n]$. Their dependence on the maximal sample size $n$
is studied in Figures 4 and 5. 
Controlling the rate of false alarms is analogous to the probability of Type I error 
for the problems of change-point detection and testing. 

Our results apply to finite samples $Y_1,\ldots,Y_n$ or sequences with a finite horizon.
For the infinite-horizon problems, thresholds
can be derived to control the average run length $\ARL$ that represents the mean time 
to a false alarm $\MTBFA=\ARL(F)=\E_F T_h$ and the mean detection delay $\MD=\ARL(G)=\E_G(T_h-\nu)^+$, 
which correspond to the Type I and Type II errors, respectively. Determining such thresholds
is called a "converse problem" in \cite{Lotov92}.  

Inequalities and asymptotic expressions (as $h\to\infty$) for $\MTBFA(T_h)$ and $\MD(T_h)$ can be found
in \cite{Lotov92}, \cite[sec. 8.2]{TartakovskyNikiforovBasseville14}, \cite[sec. 5.2]{BassevilleNikiforov93},
and references therein.

\subsection{Applications in queuing systems}

As mentioned in the introduction, CUSUM processes are popular models for waiting times and stored commodities in queuing systems
(see \cite{Feller71v2}, sec. 6.9; \cite{Kendall53}; \cite{Lindley1952queuing}). In general G/G/1 queues with one server,
the waiting time $W_{n+1}$ of the next job or customer is represented by the previous waiting time $W_n$, plus the previous service time
$S_n$, reduced by the interarrival time $A_{n+1}$. That is unless $W_n+S_n < A_{n+1}$, in which case the previous job ends
before the new job arrives, and the new job goes directly to the server with no waiting, 
\[
W_{n+1} = \l\ba{ll} W_n + S_n - A_{n+1}, & \mbox{ if } W_n + S_n - A_{n+1} \ge 0 \\ 0, & \mbox{ otherwise}.\ea\right.
\]
This matches the definition \num{sec:introduction} of a CUSUM process, with the increments $Y_{n+1}=S_n - A_{n+1}$.

Further, the case of negative increments that we considered in this paper, $\E(Y_{n+1})<0$, implies that the expected 
interarrival time exceeds the expected service time, which makes the {\sl utilization} of the system less than 1, preventing 
this queuing system from overloading. 

The studied probabilities of deviations have the following meaning in the queuing theory. The probability $\P\l \max_{[0,n]} W_t\ge h\r$
represents the chance for the waiting time, the amount of a commodity, or the length of a queue to exceed $h$ {\sl at least once}
during the time interval of length $n$. This includes undesired events such as the number of jobs exceeding the size of a buffer, or
even catastrophic events such as the amount of water in a reservoir exceeding the maximum allowed level and flowing over levees. 

\subsection{Applications in biological sequence alignment}

Another application of the CUSUM process and a threshold for its maximum is in the evaluation of alignment of biological sequences. 

The Smith–Waterman score matrix for the sequence alignment \cite[Sec. 2.3]{Durbin98},
\[
F(i,j) = \max\left\{\begin{array}{l}
0 \\
F(i-1,j-1) + \mbox{match/mismatch score} \\
F(i-1,j) - \mbox{gap penalty} \\
F(i,j-1) - \mbox{gap penalty} \end{array}\right.
\]
actually defines a two-dimensional analogue of a classical CUSUM process. In both cases, the statistic evolves through additive updates driven by incoming observations, with a reset to zero whenever the accumulated score becomes negative. This reset mechanism ensures that only sustained positive evidence contributes to detection, whether of a distributional change in CUSUM or of a locally high-similarity region in sequence alignment. In the Smith–Waterman algorithm, the maximum value of $F(i,j)$ over all matrix entries serves as a natural detection statistic, and the threshold plays a crucial role in distinguishing statistically significant alignments from spurious matches arising by chance. Similarly to the change-point detection, where threshold selection determines the trade-off between false alarms and detection delay, the choice of alignment score threshold governs the sensitivity and specificity of detecting biologically meaningful local alignments.

\appendix
\section{An alternative proof of Spitzer formula}
\label{appA}

We present an alternative proof of (\ref{equ:spz}) by mimics arguments in \cite{Wendell58} of simple proof of Spitzer's identity. The CUSUM process can be rewritten recursively 
\be
W_{n+1}=(W_n+Y_{n+1})\vee 0.
\nnn
Hence, $\{W_j\}_{j\in\Nb}$ is a homogeneous Markov process. Following \cite{Wendell58}, introduce the Banach algebra $\MA_1$ of linear combinations of cumulative distribution functions with the distribution function $e$ of the unit mass at the origin as the identity, convolution as multiplication and norm given by the total variation. We define the following linear operators:
$Fg=f\cdot g$, where  $f$ is the distribution function of $X_1$ and $Hg(x)=g(x)\indi_{\{x\geq 0\}}$, $g\in\MA_1$.  The distribution function $g_n$ of the random variable $W_n$ can be obtained recursively
$
g_k=HF g_{k-1}
$ 
for $k\in\Nb$, and therefore, 
$
g_n=(HF)^n e 
$
for all $n\in\Nb$. Then
\begin{equation}
\label{equ:geom}
\sum_{n=0}^{\infty} t^n g_n=\sum_{n=0}^{\infty} t^n (HF)^n e=(I-tHF)^{-1}e,
\end{equation}
for $t:|t|<1$, where $I$ is the identity. Note that $H$ is a projector ($H^2=H$) which projecting mass of all negative values to zero. Hence, for any bounded Borel function the following relation 	holds:
\begin{equation}
\label{equ:trint}
\int h(x)dHg(x)= \int h(x\vee 0)dg(x).
\end{equation}
Set $\tilde g=(I-tHF)^{-1}e$. Then
\begin{equation}
\label{equ:c1}
H\tilde g=\tilde g
\end{equation}
since $H$ is a projector, and, therefore, 
\begin{equation}
\label{equ:c2}
e=(I-tHF)\tilde g=H(I-tF)\tilde g=H((e-tf)\cdot\tilde g).
\end{equation}
Moreover, (\ref{equ:c1}) and  (\ref{equ:c2}) characterize $\tilde g$. Let
\be
\tilde g=\exp\bigl(-H\log(e-tf)\bigr),
\nnn 
where $\exp(g)$ and $\log(e-g)$ are both determined by Taylor--Maclaurin series for  $g$: $\|g\|<1$. Since $\MA^{\circ}_1=(I-H)\MA_1$ is a closed subalgebra of $\MA_1$,
\be
(e-tf)\cdot\tilde g=\exp((I-H)\log(e-tf))=e+(I-H)g_1
\nnn
and, therefore, 
\be
H((e-tf)\cdot\tilde g)=H(e+(I-H)g_1)=e. 
\nnn
The equation (\ref{equ:c2}) holds. Substituting $\tilde g$ to (\ref{equ:geom}) we get that
\be
\sum_{n=0}^{\infty} t^n g_n=\exp\Bigl(\sum_{k=1}^{\infty} \frac{t^k}{k} H(f^k)\Bigr).
\nnn
Let
\be
\TA(g)(t)=\int e^{itx}dg(x)
\nnn
be the Fourier--Stieltjes transformation of the distribution function $g\in\MA_1$. Since $\TA$ is linear and multiplicative $\TA(g_1\cdot g_2)=\TA(g_1)\TA(g_2)$ we conclude that
\be
\sum_{n=0}^{\infty} t^n \TA(g_n)=\exp\Bigl(\sum_{k=1}^{\infty} \frac{t^k}{k} \TA\bigl(H(f^k)\bigr)\Bigr).
\nnn
From (\ref{equ:trint}) we conclude that
\be
\TA\bigl(H(f^k)\bigr)=\int e^{iux}dHf^n(x)=\int e^{iux_+}df^n(x)=\E e^{iuS_n^+}.
\nnn
Hence,
\be
\sum_{n=0}^{\infty} t^n \E e^{iuW^*_n}=\exp\Bigl(\sum_{k=1}^{\infty} \frac{t^k}{k}\E e^{iuS^+_n}\Bigr)
\nnn
and, therefore, (\ref{equ:spz}) holds.

\section{Wendell's framework for the maximum of CUSUM  fails}
\label{appB}

In order to obtain  Spitzer's formula for the running maxima of the CUSUM process, we try the same framework as in \cite{Wendell58}. 

{\it The origin.} Let $X_1,X_2,\ldots$ be independent identically distributed (i.i.d.) random variables; $X_i\sim F$, $i\in\Nb$. Introduce the following notations:
\be
S_n=\sum_{j=1}^n X_i; \quad S_n^{\pm}=(|S_n|\pm S_n)/2; \quad R_n=\max_{0\leq k\leq n} S_n; \quad  \Ybi_n=(R_n,S_n). 
\nnn
Then $ \Ybi_1, \Ybi_2,\ldots$ is a Markov process, and the following Spitzer's identity holds
\be
\sum_{n=0}^{\infty} t^n \phi_n(u+v,-v)=\exp\Bigl[ \sum_{n=0}^{\infty}\frac{t^n}{n}(\psi^{+}_n(u)+\psi^{-}_n(v)-1)\Bigr],
\nnn
and, in the particular case of $v=0$, 
\be
\sum_{n=0}^{\infty} t^n \phi_n(u) =\exp\Bigl[ \sum_{n=0}^{\infty}\frac{t^n}{n}\psi^{+}_n(u)\Bigr],
\nnn
where $\phi_n(u',v')=\E e^{i(u' R_n+v' S_n)}$, $\psi^{+}_n(u)=\E e^{iuS_n^+}$, $\psi^{-}_n(v)=\E e^{ivS_n^-}$, $\phi_n(u)=\E e^{iu' R_n}$  are the characteristic functions. 
 
{\it Wendell's framework for CUSUM maximum.}  Let 
\be
M_n=\max_{[0,n]} W_k,
\nnn
where {$\{W_k\}_{0}^{\infty}$} is the CUSUM process. Introduce another process (which we call CUSUM as well) by the origin and the following recursive relations
\be
W^*_0=0; \quad W^*_n=(W^*_n\vee 0)+X_n. 
\nnn
Note that, {$\{W^*_j\}_{0}^{\infty}$} is a homogeneous Markov process, and
\be
M_n=\max_{[0,n]} W^*_k.
\nnn
The main target is the distribution of $M_n$. 

Since the Markov property fails for the process  {$\{M_n\}_{0}^{\infty}$}, it should be easier to investigate the bivariate homogeneous Markov random process {$\{(M_n,W_n)\}_{0}^{\infty}$}. Let $\MA$ be the Banach algebra, which contains all linear combinations of two dimensional (right continuous) distribution functions, with convolution as a multiplication, norm given by the total variation and the delta measure at the origin $\delta_{0}$ as unit element. Moreover, introduce the following linear operators: $Fg=f\cdot g$, where $f$ is the distribution function of $(0,X_1)$; $Pg(x,y)=g(x,x\wedge y)$; $Hg(x,y)=g(x,y\vee 0)$, $g\in\MA$. Then, the distribution function $g_{n+1}$ of $(M_{n+1},W_{n+1})$ can be obtained recursively 
$
g_k=PFH g_{k-1}
$ 
for $k\in\Nb$, and, therefore, 
$
g_n=(PFH)^n e 
$
for $n\in\Nb$. Then
\be
\sum_{n=1}^{\infty} t^n g_n=\sum_{n=1}^{\infty} t^n (PFH)^n e=(I-tPFH)^{-1}e,
\nnn
where $I$ is the identity map. On the other hand,  the for joint distribution function $g_k^*$ of $(M_n,W_n^*)$ the following recursion holds $g^*_k=HPF g^*_{k-1}$, $k\in\Nb$. Hence, $g^*_n=(HPF)^n e$ for all $n\in\Nb$, and
\be
\sum_{n=1}^{\infty} t^n g^*_n=\sum_{n=1}^{\infty} t^n (HPF)^n e=(I-tHPF)^{-1}e. 
\nnn
The open problem is to find $g=(I-tHPF)^{-1}e$. 

Note that  $P$ projects all mass to the left of diagonal $x=y$ horizontally to the diagonal leaving mass to the right, and $H$  projects all mass of the lower half-plane $y<0$ orthogonally onto the $x$-axis. Wendel's proof of Spitzer's identity and the current case can be noted: the only difference is using here the projector $HP$ instead of $P$ in \cite{Wendell58}. The key point in the Wendell's proof is that $P\MA$ and $(I-P)\MA$ are both closed subalgebras with respect to the multiplication since, that allow to obtain the following representation 
\be
\tilde g=(I-tPF)^{-1}e=\exp(-P\log(e-f)),
\nnn
where  $\exp(g)$ and $\log(e-g)$ are defined both by the Maclaurin series with $f^k=f\cdot\ldots\cdot f$. The analogous representation
\be
\tilde g^*=(I-tHF)^{-1}e=\exp(-H\log(e-f)),
\nnn
holds (in one dimensional case) in the proof in Appendix A. Unfortunately,  $(I-HP)\MA$ consists of distributions concentrated on the set $\Rb^2\backslash\{(x,y)\in\Rb^2: y\in [0,x]\}$, and it is not closed with respect to the convolution.




 \end{document}